\documentclass[a4paper,oneside,10pt]{article}%
\usepackage{amsmath}
\usepackage{amsfonts}
\usepackage{amssymb}
\usepackage{graphicx}
\usepackage{color}
\usepackage{authblk}
\usepackage[square,numbers,sort&compress]{natbib}%
\providecommand{\U}[1]{\protect\rule{.1in}{.1in}}

\newtheorem{theorem}{Theorem}[section]

\newtheorem{corollary}[theorem]{Corollary}

\newtheorem{definition}[theorem]{Definition}

\newtheorem{example}[theorem]{Example}

\newtheorem{lemma}[theorem]{Lemma}

\newtheorem{proposition}[theorem]{Proposition}
\newtheorem{remark}[theorem]{Remark}

\newenvironment{proof}[1][Proof]{\noindent \textbf{#1.} }{\  $\Box$}
\numberwithin{equation}{section}

\begin{document}

\title{Stochastic heat equations driven by  space-time $G$-white noise under sublinear expectation}

\author[a, b]{Xiaojun Ji}

\author[b, c, d]{Shige Peng}

\affil[a]{Research Center for Mathematics and Interdisciplinary Sciences, Shandong University, Qingdao, Shandong 266237, P. R. China.}
\affil[b]{Frontiers Science Center for Nonlinear Expectations (Ministry of Education), Shandong University, Qingdao, Shandong 266237, P. R. China.}
\affil[c]{School of Mathematics, Shandong University, Jinan, Shandong
250100, P. R. China.}
\affil[d]{Zhongtai Securities Institute for Financial Studies, Shandong University, Jinan, Shandong
250100, P. R. China.}

\affil[ ]{\textit {Email: jixj@sdu.edu.cn, peng@sdu.edu.cn}}

\renewcommand*{\Affilfont}{\small\it}

\date{}

\maketitle

{
\begin{abstract}
In this paper, we study the stochastic heat equation driven by a multiplicative space-time $G$-white noise within the framework of sublinear expectations. The existence and uniqueness of the mild solution  are proved.  By  generalizing the  stochastic Fubini theorem under  sublinear expectations, we demonstrate that the mild solution  also qualifies as a weak solution. Additionally, we derive moment estimates for the solutions.
\medskip
\medskip

\noindent\textbf{Key words}. Sublinear expectation, Space-time $G$-white noise, Stochastic heat equations, Stochastic Fubini theorem

\medskip
\medskip

\noindent{\textbf{MSC(2020)}. }  
60G65, 60H15,  60H40

\end{abstract}}

\section{Introduction}
For nearly four decades, the study of stochastic partial differential equations (SPDEs) has been an active area in mathematics with wide-ranging applications in the fields of theoretical physics, fluid mechanics, mathematical finance and control theory. 
In most cases, an SPDE is a partial differential equation  perturbed by adding a random external forcing. One important example is the stochastic heat equation 
\begin{equation}\label{intro-1}
\frac{\partial u}{\partial t}=\frac{\partial^2 u}{\partial x^2}+f(u)+\sigma(u)\xi, \ t>0,x\in\mathbb{R},
\end{equation}
where $\xi$ represents the random external force and $f$ and $\sigma$ are smooth functions.
This classical equation can describe a variety of phenomena, such as the evolution of the density of a space-time population \cite{Daw72}, the propagation of electric potentials  in neurons \cite{Walsh},  the heat density in a random medium \cite{Kh14}, 
and the random motion of a polymer chain  in a liquid \cite{Hai}. In these models,  the random external force  $\xi$ is always assumed to be  a space-time white noise, which is a Gaussian process that is independent at different time and space position. This assumption is  an idealized approximation to very short-range-dependent physical noises \cite{Kh14}.

Due to its significant applications, the stochastic heat equation has been extensively studied in the literature,  either by generalizing It\^o's theory \cite{Walsh, Daw72, DaP-Z, Par75}, or using pathwise approaches \cite{CFO11, Hai}. In the 1970's,  
Pardoux \cite{Par72, Par75}   proved the well-posedness of the solution  by extending  the variational approach for deterministic partial differential equations. This variational approach has since become a general method for studying SPDEs \cite{KR77, LR15}. Concurrently,  Dawson \cite{Daw72} demonstrated the existence and continuity of the solution to stochastic heat equations with  Lipschitz continuous coefficients by studying the properties of semigroup operators. 

Since then, numerous probabilists have advanced the  SPDE theory, leading to  the development of various approaches. In the pioneering monograph \cite{Walsh},  Walsh  established a martingale measure theory and obtained the well-posedness and regularity  of the solution to the stochastic heat equation (\ref{intro-1}) with Lipschitz coefficients. In Walsh's theory, the  space-time white noise $\xi$ is defined as a martingale measure, and the solution  is a  multi-parameter random field. Dalang \cite{Dalang}  later extended Walsh's results to nonlinear equations in high spatial dimensions and proved the necessary and sufficient conditions for  stochastic heat equations to have  random field solutions. 
Building on this, stochastic heat equations have been extensively investigated, uncovering many precise properties of the solution, 
including the  regularity \cite{Sanz, Kh14}, Feynman-Kac formulas \cite{Ber94, Hu19}, asymptotic properties \cite{Chen16, Hu19} and intermittency \cite{Ber94,  Kh14}.  These results have been applied to a wide range of fields,  such as in stochastic control theory, including filter theory \cite{LZ21}.

Another important approach to study the stochastic heat equation is the semigroup approach,  systematically developed by Da Prato and Zabczyk \cite{DaP-Z}.  In the semigroup approach, the  equation is considered as an infinite-dimensional stochastic differential equation,  with the noise in the form of  cylindrical Wiener processes,  and the solution is a function-valued  process  indexed by time $t$.  Additionally, Pardoux and Peng \cite{PP94} provided a probabilistic representation for  parabolic SPDEs by introducing backward doubly stochastic differential equations.  Recently, pathwise approaches have proven to be  effective tools for studying singular SPDEs \cite{CFO11, Hai}. More references can be found in \cite{LR15,  Hai, Hu19, Kh14}, as well as in the references cited therein.

Although various approaches have been developed to solve different problems, all of these methods  are generally applicable to the systems of our concern (e.g.,  Eq. (\ref{intro-1})). The choice of a suitable  approach  is primarily a matter of mathematical convenience. In this paper, we mainly extend Walsh's martingale measure theory to investigate the stochastic heat equations under sublinear expectation.

The analysis of classical stochastic heat equations typically takes place in the  probability space, characterized by a specific probability measure. In these systems, the  random force $\xi$ is assumed to follow a fixed distribution. However, in most real-world problems, due to  inherent complexities, the distribution of the force  fluctuates, making it challenging to identify the true distribution. In such situations, random systems exhibit probabilistic and distributional model uncertainty.  To address this higher degree of  uncertainty, it becomes essential to study stochastic heat equations  within the sublinear expectation framework,  which accounts for uncertainty in the probability distribution of the noise term.

Motivated by the model uncertainty in financial markets, Peng \cite{P2007, P2010} systemically established a nonlinear expectation theory, especially the $G$-expectation theory,  which plays a more realistic role analogous to the Wiener measure in the classical probability theory. $G$-expectation can be represented as the supremum of a family of probability measures, including mutually singular measures. Consequently,  it has become a robust theoretical tool for analyzing problems involving probability uncertainty, such as financial issues with volatility ambiguity \cite{EJ13}.

To investigate stochastic heat equations within the  sublinear expectation framework, the appropriate noise should be constructed first. Peng \cite{P2007} introduces $G$-Brownian motion, which exhibits stationary and independent increments, as a typical temporal white noise. This crucial process is then  employed to investigate $G$-OU process and the (backward) stochastic differential equations driven by $G$-Brownian motion \cite{Gao09, HJPS, HJL}.  However, in contrast to the classical Brownian motion, $G$-Brownian motion is not a $G$-Gaussian process \cite{Peng23}, which means that the finite-dimensional distributions of $G$-Brownian motion are generally  not $G$-normally distributed, owing to the asymmetry of independence under  sublinear expectations. 
As a result, $G$-Brownian motion cannot be employed to characterize the space-indexed white noises on the sublinear expectation space.
To address this issue, Ji and Peng \cite{JP} introduces the concept of $G$-Gaussian random fields, spatial $G$-white noises and space-time $G$-white noises. The related stochastic integrals with respect to  $G$-white noises are also established, laying the foundation for studying  SPDEs under sublinear expectations.

SPDE systems, perturbed by noises in both temporal and spatial dimensions, often exhibit significant  model uncertainty. This necessitates a systematic study of partial differential equations driven by $G$-white noises under sublinear expectations to quantify and control this uncertainty. The main objective of this paper is to establish stochastic heat equations within the  sublinear expectation framework, particularly focusing on stochastic heat equations  driven by the multiplicative space-time $G$-white noise ($G$-stochastic heat equations). Specifically, we consider the following equation on a finite interval:
\begin{equation}
  \left\{\begin{array}{ll}
    & \frac{\partial}{\partial t} u (t, x) = \frac{\partial^2}{\partial x^2}
    u (t, x) +b(u(t,x))+ a (u (t, x)) \dot{\mathbf{W}} (t, x),\ 0 < t \leq T, 0 \leq x
    \leq L,\\
    & \frac{\partial}{\partial x} u (t, 0) = \frac{\partial}{\partial x} u
    (t, L) = 0,\ 0 < t \leq T,\\
    & u (0, x) = u_0 (x),\ 0 \leq x \leq L,
  \end{array}\right. \label{intro-nonlinear eq}
\end{equation}
where  the initial function $u_0(x)$ is  bounded and Borel measurable, the coefficients $a(x)$ and $b(x)$ are Lipschitz continuous functions, and  $\dot{\mathbf{W}} (t, x)$ stands for the generalized  derivative of the space-time $G$-white noise. 

By extending Walsh's worthy martingale measure theory in \cite{Walsh}, we first establish the well-definedness of the $G$-stochastic heat equation (\ref{intro-nonlinear eq}) and introduce definitions of  the  mild solution and weak solution  within the appropriate  random field space.  Since the sublinear expectation represents the supremum of a family of mutually singular measures, the integrable random fields must satisfy much stricter integerability and continuity conditions. Consequently, proving the existence of a random field solution  becomes more challenging. To overcome this obstacle, we  extend  the stochastic integrals with respect to  space-time $G$-white noises and demonstrate  the regularity of the stochastic integrals of Green's functions for heat equations. Subsequently, the existence and uniqueness of the mild solution of  (\ref{intro-nonlinear eq}) are proved through the Picard iteration method.  Furthermore, the moment estimates are  provided for the  solution.
Additionally, by establishing  the stochastic Fubini theorem for  stochastic integrals with respect to space-time $G$-white noises, we show that the mild solution implies the  weak solution of  (\ref{intro-nonlinear eq}).

In particular,  we also consider the linear stochastic heat equation driven by an additive space-time $G$-white noise. An explicit weak solution (or mild solution) is given, and the existence and uniqueness of the solution are presented. We believe that the establishment of $G$-stochastic heat equations in this paper will serve to random systems   involving probability uncertainty and lay a foundation for investigating more complex SPDEs within the sublinear expectation framework. At the end of the paper, we present several examples to illustrate the potential applications of $G$-stochastic heat equations.

The remainder of this paper is organized as follows. In Section 2, the sublinear expectation theory and space-time $G$-white noises are briefly reviewed. The stochastic integrals for space-time $G$-white noises are further extended in Section 3. Section 4 is devoted to the stochastic Fubini theorem  under sublinear expectation. In Section 5, we focus on the establishment of  the linear and nonlinear  $G$-stochastic heat equations, wherein we prove the existence,  uniqueness and moment estimates of the mild solution, which is also a weak solution. Finally, in Section 6,  we briefly illustrate the potential applications of  $G$-stochastic heat equations through several examples.

\section{Preliminaries}

In this section, we recall some basic notions and properties in the sublinear
expectation theory. More details can be found in {\cite{DHP, JP, P2007,   P2010}} and references therein.

\subsection{Basic notions of sublinear expectations}

Let $\Omega$ be a given nonempty set and $\mathcal{H}$ be a linear space of
real-valued functions on $\Omega$ such that if $X_1$, $\cdots$, $X_d \in
\mathcal{H}$, then $\varphi (X_1, X_2, \cdots, X_d) \in \mathcal{H}$ for each
$\varphi \in C_{l.Lip} (\mathbb{R}^d)$, where $C_{l.Lip} (\mathbb{R}^d)$
denotes the linear space of functions satisfying that for any $x, y \in
\mathbb{R}^d$,
\begin{equation*}
 | \varphi (x) - \varphi (y) | \leq C_{\varphi}  (1 + |x|^m + |y|^m)  |x -
   y|, \text{for some } C_{\varphi} > 0, m \in \mathbb{N} \text{ depending on }
   \varphi .
\end{equation*}
$\mathcal{H}$ is considered as the space of random variables. For each $X_i \in \mathcal{H}$, $1 \leq i \leq d$, $X = (X_1,
\cdots, X_d)$ is called a
$d$-dimensional random vector, denoted by $X \in \mathcal{H}^d$.

\begin{definition}
  \label{sublinear expectation}A {sublinear expectation}
  $\hat{\mathbb{E}}$ on $\mathcal{H}$ is a functional $\hat{\mathbb{E}} :
  \mathcal{H} \rightarrow \mathbb{R}$ satisfying the following properties: for
  each $X, Y \in \mathcal{H}$,
  \begin{itemize}
    \item[(i)] {Monotonicity:}\quad$\hat{\mathbb{E}} [X] \geq
    \hat{\mathbb{E}} [Y]  \text{ if } X \geq Y$;
    
    \item[(ii)] {Constant preserving:}\quad$\hat{\mathbb{E}} [c] = c
    \text{ for } c \in \mathbb{R}$; 
    
    \item[(iii)] {Sub-additivity:}\quad$\hat{\mathbb{E}}  [X + Y]
    \leq \hat{\mathbb{E}} [X] + \hat{\mathbb{E}} [Y]$;
    
    \item[(iv)] {Positive homogeneity:}\quad$\hat{\mathbb{E}}
    [\lambda X] = \lambda \hat{\mathbb{E}} [X]  \text{ for } \lambda \geq 0$.
  \end{itemize}
  The triplet $(\Omega, \mathcal{H}, \hat{\mathbb{E}})$ is called a
  {sublinear expectation space}. In particular, if only the properties (i) and (ii) are satisfied,
  then $\hat{\mathbb{E}}$ is called a {nonlinear expectation} and  $(\Omega, \mathcal{H}, \hat{\mathbb{E}})$ is called a
  {nonlinear expectation space}.
\end{definition}

In this paper, we are mainly interested in sublinear expectations satisfying
the following regular condition:
\begin{itemize}
  \item[(v)] {Regularity:}
  {\hspace{0.17em}}{\hspace{0.17em}}{\hspace{0.17em}} If $\{X_i \}_{i =
  1}^{\infty} \subset \mathcal{H}$ satisfies that $X_i (\omega) \downarrow 0$
  as $i \rightarrow \infty$, for each $\omega \in \Omega$, then
  \begin{equation*}
   \lim_{i \rightarrow \infty}  \hat{\mathbb{E}} [X_i] = 0. 
  \end{equation*}
\end{itemize}

Within the framework of  sublinear expectations,  important definitions, including the distribution and independence of random variables,  are introduced as follows.
\begin{definition}
 For each given $d$-dimensional random vector $X$ on the sublinear
  expectation space $(\Omega, \mathcal{H}, \hat{\mathbb{E}})$, define a
functional  $\mathbb{F}_X$ on $C_{l.Lip} (\mathbb{R}^d)$ by
\begin{equation*} \mathbb{F}_X [\varphi] := \hat{\mathbb{E}} [\varphi (X)], \text{ for
   } \varphi \in C_{l.Lip} (\mathbb{R}^d) .
\end{equation*}
$\mathbb{F}_X$ is called
the {distribution} of $X$.
\end{definition}
\begin{definition}
  Two $d$-dimensional random vectors $X_1$ and $X_2$,  defined on sublinear
  expectation spaces $(\Omega_1, \mathcal{H}_1, \hat{\mathbb{E}}_1)$ and
  $(\Omega_2, \mathcal{H}_2, \hat{\mathbb{E}} \mathbb{}_2)$, respectively, are
  called {identically distributed}, denoted by $X_1 \overset{d}{=}
  X_2$, if $\mathbb{F}_{X_1} =\mathbb{F}_{X_2}$, i.e.,
  \begin{equation*} \hat{\mathbb{E}}_1 [\varphi (X_1)] = \hat{\mathbb{E}}_2 [\varphi (X_2)], 
     \text{ for } \varphi \in C_{l.Lip} (\mathbb{R}^d) .
  \end{equation*}
\end{definition}

\begin{definition}
 An $n$-dimensional random vector $Y$ is said to be
  {independent} of another $m$-dimensional random vector $X$ on $(\Omega, \mathcal{H}, \hat{\mathbb{E}})$ if, for each test function $\varphi \in
  C_{l.Lip} (\mathbb{R}^{m + n})$,
  \begin{equation*} 
  \hat{\mathbb{E}} [\varphi (X, Y)] = \hat{\mathbb{E}} [\hat{\mathbb{E}}
     [\varphi (x, Y)]_{x = X}] . 
  \end{equation*}
  In particular, let $\bar{X}$ and $X$ be two $m$-dimensional random vectors on $(\Omega,
  \mathcal{H}, \hat{\mathbb{E}})$. $\bar{X}$ is called an {independent copy} of
  $X$ if $\bar{X} \overset{d}{=} X$ and $\bar{X}$ is independent of $X$.
\end{definition}

\begin{remark}
  It is important to mention that ``$Y$ is independent of $X$'' does not imply
  that ``$X$ is independent of $Y$'' (See  Peng {\cite{P2010}}).
\end{remark}

\subsection{$G$-white noises under sublinear expectations}
Within the  sublinear expectation framework, essential random fields,  including  $G$-Brownian motion and $G$-white noises, are introduced in \cite{P2007, JP, P2010}.

According to Peng's central limit theorem under sublinear expectations \cite{P2010, K20, Song20}, in the world of probability model uncertainty, $G$-normal distribution plays the same central role as the normal distribution in the classical cases.

\begin{definition}
  A $d$-dimensional random vector $X$ on a sublinear
  expectation space $(\Omega, \mathcal{H}, \hat{\mathbb{E}})$ is called
  (centralized) {$G$-normally distributed} if
  \begin{equation*}
   aX + b \bar{X} \overset{d}{=} \sqrt{a^2 + b^2} X, \text{ for } a, b \geq 0,
  \end{equation*}
  where $\bar{X}$ is an independent copy of $X$.
\end{definition}

\begin{definition}
  On a sublinear expectation space
  $(\Omega, \mathcal{H}, \hat{\mathbb{E}})$, a $d$-dimensional stochastic process $(B_t)_{t \geq 0}$ with $B_t\in \mathcal{H}^d$ for $t\geq0$  is called a {(symmetric) $G$-Brownian motion} if
  the following properties are satisfied:
  \begin{itemize}
    \item[(i)] $B_0 = 0 ;$
    
    \item[(ii)] For each $t, s \geq 0,$ $B_{t + s} - B_t$ and $B_s$ are
    identically distributed and $B_{t + s} - B_t$ is independent of $(B_{t_1}, \cdots, B_{t_n})$ for any $n \in \mathbb{N}$ and $0 \leq t_1 \leq \cdots \leq t_n \leq t.$
    
    \item[(iii)] $\hat{\mathbb{E}} [B_t] = \hat{\mathbb{E}} [- B_t] = 0$ and $\lim_{t
    \downarrow 0} \hat{\mathbb{E}} [| B_t |^3] t^{- 1} = 0$.
  \end{itemize}
\end{definition}

\begin{definition}
  Let $\Gamma$ be any parameter set. A family of random variables $(W_{\gamma})_{\gamma \in \Gamma}$ on $(\Omega, \mathcal{H}, \hat{\mathbb{E}})$ is
  called a $1$-dimensional {$G$-Gaussian random field} if for each $n\in\mathbb{N}$ and
  $\gamma_1, \cdots, \gamma_n \in
  \Gamma$,  $n$-dimensional random vector
  $ (W_{\gamma_1}, \cdots, W_{\gamma_n}) \in
  \mathcal{H}^n$ is $G$-normally distributed.
\end{definition}

 It is worth noting that for any given $t>0$,  $B_t$ is $G$-normally distributed. However, the finite-dimensional distributions of $(B_t)_{t\geq0}$ are generally not $G$-normally distributed, implying that $G$-Brownian motion is no longer a  $G$-Gaussian random field.

Let $0 \leq \underline{\sigma}^2 \leq \overline{\sigma}^2$ be any given numbers, and  $\mathcal{B}_0 (\mathbb{R}^d) := \{A
  \in \mathcal{B}(\mathbb{R}^d) : \lambda_A < \infty\}$, where $\lambda_A$
  denotes the Lebesgue measure of $A \in \mathcal{B} (\mathbb{R}^d)$.  

\begin{definition}
  \label{Gwhitenoise}A
  $1$-dimensional  $G$-Gaussian random field $\mathbb{W} = (\mathbb{W}_A)_{A
  \in \mathcal{B}_0 (\mathbb{R}^d)}$ on $(\Omega, \mathcal{H}, \hat{\mathbb{E}})$ is called a {spatial $G$-white
    noise} if
  \begin{itemize}
    \item[(i)] For each $A \in \mathcal{B}_0 (\mathbb{R}^d)$,
    $\hat{\mathbb{E}} [\mathbb{W}_A^2] = \overline{\sigma}^2 \lambda_A$, $-
    \hat{\mathbb{E}} [- \mathbb{W}_A^2] = \underline{\sigma}^2 \lambda_A$;
    
    \item[(ii)] For each $A_1, A_2 \in \mathcal{B}_0 (\mathbb{R}^d)$ such that $A_1
    \cap A_2 = \emptyset$,
    \begin{align*}
      & \hat{\mathbb{E}} [\mathbb{W}_{A_1} \mathbb{W}_{A_2}] =
      \hat{\mathbb{E}} [- \mathbb{W}_{A_1} \mathbb{W}_{A_2}] = 0, \\
      & \hat{\mathbb{E}} [(\mathbb{W}_{A_1 \cup A_2} - \mathbb{W}_{A_1} -
      \mathbb{W}_{A_2})^2] = 0. 
    \end{align*}
  \end{itemize}
\end{definition}

In the following, we consistently fix the parameters  $\underline{\sigma}^2$ and   $\overline{\sigma}^2$. 
It is important to note that  although the distribution of the spatial $G$-white noise has only two parameters $\overline{\sigma}^2$ and $\underline{\sigma}^2$, it is not uniquely determined by these parameters. According to Theorem 4.4 in Ji and Peng \cite{JP}, its distribution is characterized  by a family of sublinear generating functions $\{G_{A_1,\dots,A_n}(\cdot):\forall n\in\mathbb{N}, A_1,\dots,A_n\in \mathcal{B}_0(\mathbb{R}^d)\}$. Different choices of this family
of generating functions can lead to distinct  spatial $G$-white
noises with the same parameters $\overline{\sigma}^2$ and $\underline{\sigma}^2$.


Building upon the $G$-Brownian motion and  spatial $G$-white noise, we introduce a  spatio-temporal random field, which is called the space-time $G$-white noise.
Set the parameter set as
\begin{equation*}
 \Gamma = \{ [s, t) \times A : 0 \leq s \leq t < \infty, A \in \mathcal{B}_0
   (\mathbb{R}^d) \} . 
\end{equation*}
\begin{definition}
  \label{tx-Gwhitenoise}A random field $\{ \mathbf{W} ([s, t) \times
  A)\}_{([s, t) \times A) \in \Gamma}$ on a sublinear expectation space
  $(\Omega, \mathcal{H}, \hat{\mathbb{E}})$ is called a $1$-dimensional
  {space-time $G$-white noise} if it satisfies the following conditions:
  \begin{itemize}
    \item[(i)] For each fixed $[s, t)$, the random field $\{ \mathbf{W} ([s,
    t) \times A)\}_{A \in \mathcal{B}_0 (\mathbb{R}^d)}$ is a $1$-dimensional
    spatial $G$-white noise and has the same family of finite-dimensional
    distributions as $(\sqrt{t - s} \mathbb{W}_A)_{A \in \mathcal{B}_0
    (\mathbb{R}^d)}$;
    
    \item[(ii)] For any $r \leq s \leq t$, $A \in \mathcal{B}_0
    (\mathbb{R}^d)$, $\mathbf{W} ([r, s) \times A) + \mathbf{W} ([s, t) \times
    A) = \mathbf{W} ([r, t) \times A)$;
    
    \item[(iii)] For any $0\leq s_i \leq t_i \leq s\leq t$ and $A, A_i \in \mathcal{B}_0 (\mathbb{R}^d)$,  $i = 1,
    \cdots, n$, $n\in\mathbb{N}$, $\mathbf{W} ([s, t) \times A)$ is independent of $(\mathbf{W}
    ([s_1, t_1) \times A_1), \cdots, \mathbf{W} ([s_n, t_n) \times A_n))$,    
  \end{itemize}
  where $(\mathbb{W}_A)_{A \in \mathcal{B}_0 (\mathbb{R}^d)}$ is a
  $1$-dimensional spatial $G$-white noise with parameters $\underline{\sigma}^2$ and $\overline{\sigma}^2$.
\end{definition}

Denote 
\begin{align*}
  \Gamma_1 = &\{ [s, t) \times (x,y): (x,y)=(x_1,y_1)\times\cdots\times(x_d,y_d), \forall0 \leq s \leq t < \infty,  x=(x_1,\cdots,x_d), \\
   &\ y=(y_1,\cdots,y_d)\in 
   \mathbb{R}^d\text{ such that } -\infty<x_i\leq y_i<\infty \text{ for }1\leq i\leq d\}.
\end{align*}
Then $\Gamma_1\subset\Gamma$. We first construct a space-time $G$-white noise indexed by $\Gamma_1$ and then extend it to the index set $\Gamma$.
Let
\begin{equation*} 
\Omega = \{ \omega \in C_0(\Gamma_1;\mathbb{R}): \omega (A_1\cup A_2) =
   \omega (A_1) + \omega (A_2),\text{ for } A_1, A_2\in\Gamma_1 \text{ with } A_1\cap A_2=\emptyset\}
\end{equation*}
 be the space of $\mathbb{R}$-valued continuous functions  endowed with the distance
  \begin{equation*}
   d(\omega,\tilde{\omega})=\sum_{k=1}^\infty2^{-k}\left[\left(\max_{\substack{0\leq s\leq t\leq k\\ (x,y)\subset[-k,k]^d}}|\omega([s,t)\times(x,y))-\tilde{\omega}([s,t)\times(x,y))|\right)\wedge1\right],\forall\omega,\tilde{\omega}\in\Omega.
   \end{equation*}

For each $\omega \in \Omega$, define the canonical process
$(\mathbf{W}_{\gamma})_{\gamma \in \Gamma_1}$ by
\begin{equation*} 
\mathbf{W} ([s, t) \times (x,y)) (\omega) = \omega ([s, t) \times (x,y)), \ \forall\, 0
   \leq s \leq t < \infty, (x,y) \in \mathcal{B}_0 (\mathbb{R}^d) . 
\end{equation*}
For any $T \geq 0$, set $\mathcal{F}_T = \sigma \{ \mathbf{W} ([s, t) \times
(x,y)), 0 \leq s \leq t \leq T, (x,y)\in \mathcal{B}_0 (\mathbb{R}^d)\}$,
$\mathcal{F}= \bigvee_{T \geq 0} \mathcal{F}_T$, and
\begin{align*}
  L_{ip} (\mathcal{F}_T) := & \{\varphi (\mathbf{W} ([s_1, t_1) \times (x_1,y_1)),
  \cdots, \mathbf{W} ([s_n, t_n) \times (x_n,y_n))) : \text{for any } n \in
  \mathbb{N}, \\
  &\hspace{0.17em} \hspace{0.17em} \hspace{0.17em} \varphi \in C_{l.Lip}(\mathbb{R}^n), 0\leq s_i \leq t_i \leq T,(x_i,y_i) \in \mathcal{B}_0 (\mathbb{R}^d) \text{ for } i = 1, \cdots, n  \},\\
  L_{ip} (\mathcal{F}) := & \cup_{n = 1}^{\infty} L_{ip} (\mathcal{F}_n) .
\end{align*}

In {\cite{JP}}, we construct a sublinear expectation $\hat{\mathbb{E}}$ and a conditional expectation $\hat{\mathbb{E}}[\cdot|\mathcal{F}_t]$ on
$L_{ip} (\mathcal{F})$, under which the canonical process $(\mathbf{W}_\gamma)_{\gamma\in\Gamma_1}$ is a
$1$-dimensional space-time $G$-white noise. For each $p \geq 1$, $T \geq
0$, denote by $\mathbf{L}_G^p (\mathcal{F}_T)$ (resp., $\mathbf{L}_G^p
(\mathcal{F})$) the completion of $L_{ip} (\mathcal{F}_T)$ (resp., $L_{ip}
(\mathcal{F})$) under the norm $\|\cdot\|_p := (\hat{\mathbb{E}} [| \cdot
|^p])^{1 / p}$. Then the sublinear expectation and conditional expectation can be extended to the complete space $\mathbf{L}_G^p(\mathcal{F})$.

For any $\gamma=[s,t)\times A\in\Gamma$,  there exists a sequence $\{B_m\}_{m=1}^\infty$ of finite unions of open intervals such that $\lim_{m\rightarrow \infty }\lambda_{A\triangle B_m}=0$.
By Definition \ref{tx-Gwhitenoise}, we know that $\{\mathbf{W}([s,t)\times B_m)\}_{m=1}^\infty$ is a Cauchy sequence in $\mathbf{L}_G^p(\mathcal{F}_t)$ and denote the limit by $\mathbf{W}([s,t)\times A)\in \mathbf{L}_G^p(\mathcal{F}_t)$. Thus, we construct a spatio-temporal random field $\mathbf{W}=(\mathbf{W}_\gamma)_{\gamma\in\Gamma}$ indexed by $\Gamma$  and it is easy to verify that $\mathbf{W}$  is a space-time $G$-white noise  under $\hat{\mathbb{E}}$. 

Furthermore, denote
\begin{align*}
  \mathcal{H}_T := & \{\varphi (\mathbf{W} ([s_1, t_1) \times A_1),
  \cdots, \mathbf{W} ([s_n, t_n) \times A_n)) : \text{for any } n \in
  \mathbb{N}, \hspace{0.17em} \varphi \in C_{l.Lip}(\mathbb{R}^n), \\
  & \hspace{0.17em} \hspace{0.17em} \hspace{0.17em} 0\leq s_i \leq t_i \leq T \text{ for }i = 1, \cdots, n,A_1,
  \cdots, A_n \in \mathcal{B}_0 (\mathbb{R}^d)  \},\\
 \mathcal{H} := & \cup_{n = 1}^{\infty} \mathcal{H}_n .
\end{align*}
Then $L_{ip}
(\mathcal{F}_T)\subset\mathcal{H}_T\subset\mathbf{L}_G^p(\mathcal{F}_T)$ and $L_{ip}(\mathcal{F})\subset\mathcal{H}\subset\mathbf{L}_G^p(\mathcal{F})$. The completion of $\mathcal{H}_T$ (resp., $\mathcal{H}$) under the norm $\|\cdot\|_p$ is equal to $\mathbf{L}_G^p(\mathcal{F}_T)$ (resp., $\mathbf{L}_G^p
(\mathcal{F})$).

It is worth mentioning that  random variables in the Banach space $\mathbf{L}_G^p(\mathcal{F})$ can be treated as a quasi-continuous function defined on $\Omega$.  In fact, by generalizing Theorem 6.6 in  \cite{JP}, we have the following robust representation
theorem for  sublinear expectations.

\begin{theorem}\label{rep-thm}(See \cite{DHP, HP, JP})
There exists a weakly compact set $\mathcal{P}$ of probability measures on $(\Omega,\mathcal{B}(\Omega))$ such that
  \begin{equation*} 
  \hat{\mathbb{E}} [X] = \max_{P \in \mathcal{P}} E_P [X],  \text{ for all }
     X \in L_{ip}({\mathcal{F}}),
  \end{equation*}
 where $\mathcal{B}(\Omega)$ is the $\sigma$-algebra generated by all open sets.  $\mathcal{P}$ is called the  family of
  probability measures that represents $\hat{\mathbb{E}}$.
\end{theorem}

Let $\mathcal{P}$ be the weakly compact set that represents $\hat{\mathbb{E}}$.
For this set $\mathcal{P}$, the related capacity can be defined by
\begin{equation*} 
c (A) := \sup_{P \in \mathcal{P}} P (A), \text{ for }A \in \mathcal{B}(\Omega).
\end{equation*}
A set $A \in \mathcal{B}(\Omega)$ is called polar if $c (A) = 0$. A
property holds ``quasi-surely'' (q.s.) if it holds outside a polar set. In the
following, we do not distinguish two random variables $X$ and $Y$ if $X =Y$ q.s.

\begin{definition}
A random variable $X$ on $\Omega$ is called quasi-continuous if
for each $ \epsilon> 0$, there exists an open set $O$ with $c(O) <\epsilon$ such that $X|_{O^c}$ is continuous.

Furthermore, we say  $X$ has a quasi-continuous version if there exists a quasi-continuous
random variable $Y$ such that $ X = Y$ q.s.
\end{definition}

Let $ L^0(\Omega)$ denote the space of all $\mathcal{B}(\Omega)$-measurable real functions. We can define the related upper expectation for $X\in L^0(\Omega)$, making the following definition meaningful,
\begin{equation*}
\bar{\mathbb{E}}[X]:=\sup_{P\in\mathcal{P}}E_P[X], \text{ for } X\in L^0(\Omega).
\end{equation*}
Then $\hat{\mathbb{E}}=\bar{\mathbb{E}}$ on $L_{ip}({\mathcal{F}})$ and the completions of $L_{ip}({\mathcal{F}})$ under norms $(\hat{\mathbb{E}} [| \cdot
|^p])^{1 / p}$ and $(\bar{\mathbb{E}} [| \cdot
|^p])^{1 / p}$ are the same. By Theorem 25 in Denis et al. \cite{DHP} and Proposition 4.4 in Hu and Peng \cite{HP}, we have the following pathwise characterization of $\mathbf{L}_G^p(\mathcal{F})$,
\begin{equation*}
\mathbf{L}_G^p(\mathcal{F})=\left\{X\in L^0(\Omega): X \text{ has a quasi-continuous version and }\lim_{N\rightarrow\infty}\bar{\mathbb{E}}[|X|^p\mathbf{1}_{\{|X|>N\}}]=0 \right\}.
\end{equation*}

In particular, if $X= Y$  in $\mathbf{L}_G^p(\mathcal{F})$, then $X = Y$ q.s. For $ X\in \mathbf{L}_G^p(\mathcal{F})$ and $\{X_n\}_{n=1}^\infty\subset\mathbf{L}_G^p(\mathcal{F})$, we say that $\{X_n\}_{n=1}^\infty$ converges to $X$ in $\mathbf{L}_G^p(\mathcal{F})$, denoted by $X=\mathbb{L}^p-\lim_{n\rightarrow\infty}X_n$, if $\lim_{n\rightarrow\infty}\hat{\mathbb{E}}[|X_n-X|^p] = 0$.

\section{Stochastic integrals with respect to space-time $G$-white noises}
For each given $p\geq1$, $T>0$, and $K\in\mathcal{B}(\mathbb{R}^d)$,
let $\mathbf{M} ^{p, 0} ([0, T] \times K)$ be the collection of
simple random fields with the form:
\begin{equation}
  \eta (t, x ; \omega) := \sum_{i = 0}^{n - 1} \sum_{j = 1}^m X_{ij} (\omega) 
    \textbf{1}_{[t_i, t_{i + 1})} (t)\textbf{1}_{A_j} (x), \label{simple-func}
\end{equation}
where $X_{ij} \in \mathbf{L}_G^p (\mathcal{F}_{t_i})$, $i = 0, \cdots, n -
1$, $j = 1, \cdots, m$, $0 = t_0 < t_1 < \cdots < t_n = T$, and $\{A_j \}_{j =
1}^m \subset \mathcal{B}_0 (\mathbb{R}^d)$ is a mutually disjoint sequence such that $\cup_{j=1}^mA_j\subset K$. The Bochner integral and stochastic integral for $\eta\in\mathbf{M} ^{p, 0} ([0, T] \times K)$ can be defined as follows:
\begin{align}
\int_0^T\int_{K}\eta (t, x)dxdt&:=\sum_{i = 0}^{n - 1} \sum_{j = 1}^m X_{ij} (t_{i+1}-t_i)\lambda_{A_j},\label{bointeg}\\
\int_0^T\int_{K}\eta (t, x)\mathbf{W} (dt,
    dx)&:=\sum_{i = 0}^{n - 1} \sum_{j = 1}^m X_{ij} \mathbf{W} ([t_i,t_{i+1})\times A_j).\label{stinteg}
\end{align}
These integrals (\ref{bointeg}) and (\ref{stinteg}) are linear mappings from $\mathbf{M} ^{p, 0} ([0, T] \times K)$ to $\mathbf{L}_G^p(\mathcal{F}_T)$. 
Similarly, for any $0\leq r\leq s\leq T$ and $K'\subset K$, we can also define 
\begin{align}
\int_r^s\int_{K'}\eta (t, x)dxdt&:=\int_0^T\int_{K}\eta (t, x)\mathbf{1}_{K'}(x)\mathbf{1}_{[r,s]}(t)dxdt,\label{bointeg-2}\\
 \int_r^s\int_{K'}\eta (t, x)\mathbf{W} (dt,dx)&:=\int_0^T\int_{K}\eta (t, x)\mathbf{1}_{K'}(x)\mathbf{1}_{[r,s]}(t)\mathbf{W} (dt, dx).
\end{align}
Set
\begin{equation*} \|\eta\|_{\mathbf{M}^p} := \left( \hat{\mathbb{E}} \left[ \int_0^T
   \int_{K} |\eta (t, x) |^p dxdt \right] \right)^{\frac{1}{p}},\ 
   \|\eta\|_{\mathbf{S}^p} := \sup_{t \in [0, T]} \sup_{x \in K}
   (\hat{\mathbb{E}} [| \eta (t, x) |^p])^{\frac{1}{p}}, 
\end{equation*}
and denote by $\mathbf{M}_G^p ([0, T] \times K)$ and  $\mathbf{S}_G^p
([0, T] \times K)$ the completions of $\mathbf{M} ^{p, 0} ([0, T]
\times K)$ under the norm $\| \cdot \|_{\mathbf{M}^p}$ and $\|
\cdot \|_{\mathbf{S}^p}$, respectively. The integral (\ref{bointeg}) (resp., (\ref{stinteg})) can be continuously extended to the complete space $\mathbf{M}_G^1 ([0, T] \times K)$ (resp., $\mathbf{M}_G^2([0, T] \times K)$). In particular, if $K\in \mathcal{B}_0(\mathbb{R}^d)$, then the Bochner integral in Eq. (\ref{bointeg}) can be continuously extended to  $\mathbf{M}_G^2 ([0, T] \times K)$.

For the space-time $G$-white noise and related stochastic integrals, we have the following important and useful properties.
\begin{lemma}\label{stointprop}(See Ji and Peng \cite{JP})
For each $0\leq s\leq r\leq T$, we have the following properties:
  \begin{itemize}
    \item[(i)] For each $\xi\in \mathbf{L}_G^2 (\mathcal{F}_s)$ and $A,\bar{A}\in\mathcal{B}_0({\mathbb{R}^d})$ such that $A\cap\bar{A}=\emptyset$,
    \begin{align}
  & \hat{\mathbb{E}} \left[ \xi\mathbf{W} ([s,r)\times A) \right] = 0,\label{mean0}\\
  & \hat{\mathbb{E}} \left[ \xi\mathbf{W} ([s,r)\times A)\mathbf{W} ([s,r)\times \bar{A}) \right] = 0,\label{dismean0}
  \end{align}
  \item[(ii)] For each $\eta\in\mathbf{M}^2_G ([0, T] \times \mathbb{R}^d)$,
  \begin{align}
  & \hat{\mathbb{E}} \left[ \int_0^T \int_{\mathbb{R}^d} \eta (t, x) \mathbf{W} (dt,
    dx) \right] = 0,  \\
    & \hat{\mathbb{E}} \left[ \left|\int_0^T \int_{\mathbb{R}^d} \eta (t, x) \mathbf{W} (dt,
    dx) \right|^2\right] \leq \overline{\sigma}^2\hat{\mathbb{E}} \left[ \int_0^T \int_{\mathbb{R}^d} |\eta (t, x)|^2 dxdt \right].
  \end{align}
  \end{itemize}

\end{lemma}

\begin{lemma}
  \label{nle-1}For each $0 \leq s \leq r \leq T$, $\xi\in\mathbf{L}_G^2(\mathcal{F}_s)$, $\eta \in
  \mathbf{M}^2_G ([0, T] \times \mathbb{R}^d)$, and $K_1, K_2 \in
  \mathcal{B}_0 (\mathbb{R}^d)$ such that $K_1 \cap K_2 = \emptyset$, we have
  \begin{align}
    & \hat{\mathbb{E}} \left[ \xi \left( \int_s^r \int_{K_1} \eta (t, x) \mathbf{W} (dt, dx)
    \right) \right] = 0,  \label{f-eq2}\\
    & \hat{\mathbb{E}} \left[ \pm \left( \int_0^s \int_{K_1} \eta (t, x)
    \mathbf{W} (dt, dx) \right) \left( \int_0^s \int_{K_2} \eta (t, x) \mathbf{W}
    (dt, dx) \right) \right] = 0.  \label{f-eq3}
  \end{align}
\end{lemma}

\begin{proof}
  Equality (\ref{f-eq2}) is direct from Lemma \ref{stointprop} and Proposition 5.8 in  \cite{JP}. It suffices to prove Eq. (\ref{f-eq3}) here. We firstly
  verify that Eq. (\ref{f-eq3}) holds for $\eta \in \mathbf{M}^{2, 0}([0, T] \times
  \mathbb{R}^d)$. Without loss of generality, assume that $\eta (t, x) =
  \sum_{i = 0}^{n - 1} \sum_{j = 1}^m X_{ij} \text{\textbf{1}}_{[t_i, t_{i +
  1})} (t) \text{\textbf{1}}_{A_j} (x)$, where $0 \leq t_0 <\cdots <
  t_n \leq T$, $\{A_j\}_{j=1}^m\subset \mathcal{B}_0 (\mathbb{R}^d)$ is a mutually disjoint sequence, and $X_{ij}
  \in \mathbf{L}_G^2 (\mathcal{F}_{t_i})$, $i = 0, \cdots, n - 1$, $j = 1,
  \cdots, m$. Then
  \begin{align*}
    & \hat{\mathbb{E}} \left[ \left( \int_0^s \int_{K_1} \eta (t, x) \mathbf{W}
    (dt, dx) \right) \left( \int_0^s \int_{K_2} \eta (t, x) \mathbf{W} (dt, dx)
    \right) \right]\\
    & = \hat{\mathbb{E}} \left[ \sum_{i, l = 0}^{n - 1} \sum_{j, k = 1}^m
    X_{ij} X_{lk} \mathbf{W} ([t_i \wedge s, t_{i + 1} \wedge s) \times (A_j
    \cap K_1)) \mathbf{W} ([t_l \wedge s, t_{l + 1} \wedge s) \times (A_k \cap
    K_2)) \right] .
  \end{align*}
  By (i) in Lemma \ref{stointprop},
  \begin{align*}
    & \hat{\mathbb{E}} \left[ \left( \int_0^s \int_{K_1} \eta (t, x) \mathbf{W}
    (dt, dx) \right) \left( \int_0^s \int_{K_2} \eta (t, x) \mathbf{W} (dt, dx)
    \right) \right]\\
    & = \hat{\mathbb{E}} \left[ \sum_{i = 0}^{n - 1} \sum_{j,k = 1}^m X_{ij} X_{ik} \mathbf{W} ([t_i \wedge s, t_{i + 1} \wedge s) \times
    (A_j \cap K_1)) \mathbf{W} ([t_i \wedge s, t_{i + 1} \wedge s) \times (A_k
    \cap K_2)) \right]\\
    & = 0.
  \end{align*}
  Hence, Eq. (\ref{f-eq3}) holds for $\eta \in \mathbf{M}^{2, 0} ([0, T] \times
  \mathbb{R}^d)$.

  For $\eta \in \mathbf{M}^2_G ([0, T] \times \mathbb{R}^d)$,
  there exists a sequence of simple random fields $\{\eta_n \}_{n = 1}^{\infty} \subset \mathbf{M}^{2, 0} ([0, T]
  \times \mathbb{R}^d)$ such that $\eta_n$ converges to $\eta$ under $\| \cdot  \|_{\mathbf{M}^2}$ as $n \rightarrow \infty$. Taking $n \rightarrow \infty$,  we observe that
  \begin{align*}
    & \hat{\mathbb{E}} \left[ \left| \left( \int_0^s \int_{K_1} (\eta (t, x) -
    \eta_n (t, x)) \mathbf{W} (dt, dx) \right) \left( \int_0^s \int_{K_2} \eta (t,
    x) \mathbf{W} (dt, dx) \right) \right| \right]\\
    & \leq \hat{\mathbb{E}} \left[ \left| \int_0^s \int_{K_1} (\eta (t, x) - \eta_n
    (t, x)) \mathbf{W} (dt, dx) \right|^2 \right]^{1 / 2} \hat{\mathbb{E}}
    \left[ \left| \int_0^s \int_{K_2} \eta (t, x) \mathbf{W} (dt, dx) \right|^2
    \right]^{1 / 2}\\
    & \leq \overline{\sigma}^2 \|\eta\|_{\mathbf{M}^2}  \|\eta - \eta_n \|_{\mathbf{M}^2}
    \rightarrow 0,
  \end{align*}
  and
  \begin{align*}
     &\hat{\mathbb{E}} \left[ \left| \left( \int_0^s \int_{K_1} \eta_n (t, x)
    \mathbf{W} (dt, dx) \right) \left( \int_0^s \int_{K_2} (\eta (t, x) - \eta_n (t,
    x)) \mathbf{W} (dt, dx) \right) \right| \right]\\
    & \leq \overline{\sigma}^2 \|\eta_n \|_{\mathbf{M}^2}  \|\eta - \eta_n
    \|_{\mathbf{M}^2} \rightarrow 0.
  \end{align*}
  Hence, it follows that
  \begin{align*}
  &\hat{\mathbb{E}} \left[ \left( \int_0^s \int_{K_1} \eta (t, x) \mathbf{W}
     (dt, dx) \right) \left( \int_0^s \int_{K_2} \eta (t, x) \mathbf{W} (dt, dx)
     \right) \right] \\
     &=\lim_{n\rightarrow\infty}\hat{\mathbb{E}} \left[ \left( \int_0^s \int_{K_1} \eta_n (t, x) \mathbf{W}
     (dt, dx) \right) \left( \int_0^s \int_{K_2} \eta_n (t, x) \mathbf{W} (dt, dx)
     \right) \right]\\
     &= 0. 
  \end{align*}
  Similarly, we have
  \begin{equation*} \hat{\mathbb{E}} \left[ - \left( \int_0^s \int_{K_1} \eta (t, x) \mathbf{W}
     (dt, dx) \right) \left( \int_0^s \int_{K_2} \eta (t, x) \mathbf{W} (dt, dx)
     \right) \right] = 0, 
  \end{equation*}
  which completes the proof.
\end{proof}

Denote by $L^2  ([0, T] \times \mathbb{R}^d)$ the linear space of square-integrable functions on $[0, T] \times \mathbb{R}^d$, that is,
\begin{equation*}
L^2  ([0, T] \times \mathbb{R}^d)=\left\{h:[0, T] \times \mathbb{R}^d\rightarrow\mathbb{R} \text{ s.t. } \int_0^T\int_{\mathbb{R}^d}|h(t,x)|^2dxdt< \infty\right\}.
\end{equation*}
{It is important to note that $L^2  ([0, T] \times \mathbb{R}^d)
 \subset \mathbf{M}_G^2 ([0, T] \times \mathbb{R}^d)$.
Then for each $h (t, x) \in
L^2  ([0, T] \times \mathbb{R}^d)$ and $\eta (t, x) = \sum_{i = 0}^{n - 1} \sum_{j =
1}^m X_{ij}  \textbf{1}_{[t_i, t_{i + 1})} (t)  \textbf{1}_{A_j} (x) \in
\mathbf{M}^{2, 0} ([0, T] \times \mathbb{R}^d)$, we can define the  stochastic
integral $I_h (\eta)$ from $\mathbf{M}^{2, 0} ([0, T] \times \mathbb{R}^d)$ to
{$\mathbf{L}_G^2 (\mathcal{F}_T)$} as follows:
\begin{equation}
  I_h (\eta) = \int_0^T \int_{\mathbb{R}^d} h (t, x) \eta (t, x) \mathbf{W} (dt, dx)
:= \sum_{i = 0}^{n - 1} \sum_{j = 1}^m X_{ij}  \int_{t_i}^{t_{i + 1}}
  \int_{A_j} h (t, x) \mathbf{W} (dt, dx) . \label{stochastic integral2}
\end{equation}

For this stochastic integral $I_h$, we have analogous properties as presented in Lemma \ref{stointprop}.

\begin{lemma}
  \label{Ito ineq}For each $h \in L^2  ([0, T] \times \mathbb{R}^d), \eta
   \in \mathbf{M}^{2, 0} ([0, T] \times \mathbb{R}^d)$, we have
  \begin{align}
    \label{nneq1} &\hat{\mathbb{E}} \left[  \int_0^T \int_{\mathbb{R}^d}
    h (t, x) \eta (t, x) \mathbf{W} (dt, dx)  \right] =0,\\
    &\label{neq-1} \hat{\mathbb{E}} \left[ \left| \int_0^T \int_{\mathbb{R}^d}
    h (t, x) \eta (t, x) \mathbf{W} (dt, dx) \right|^2 \right] \leq
    \overline{\sigma}^2 || \eta ||_{\mathbf{S}^2}^2 || h ||_{L^2}^2 .
  \end{align}
\end{lemma}
{}
\begin{proof} Property (\ref{nneq1}) is obvious from the definition (\ref{stochastic integral2}) and Lemma \ref{nle-1}. We only prove Eq. (\ref{neq-1}) here.

For each $\eta \in \mathbf{M}^{2, 0} ([0, T] \times \mathbb{R}^d)$ with the
  form $\eta (t, x) = \sum_{i = 0}^{n - 1} \sum_{j = 1}^m X_{ij} 
  \textbf{1}_{[t_i, t_{i + 1})} (t)  \textbf{1}_{A_j} (x)$, by Lemma
  \ref{nle-1},
  \begin{align*}
    \hat{\mathbb{E}} \left[ \left| \int_0^T \int_{\mathbb{R}^d} h (t, x) \eta (t,
    x) \mathbf{W} (dt, dx) \right|^2 \right] & = \hat{\mathbb{E}} \left[
    \left| \sum_{i = 0}^{n - 1} \sum_{j = 1}^m X_{ij}  \int_{t_i}^{t_{i + 1}}
    \int_{A_j} h (t, x) \mathbf{W} (dt, dx) \right|^2 \right]\\
    & = \hat{\mathbb{E}} \left[ \sum_{i = 0}^{n - 1} \sum_{j = 1}^m X_{ij}^2
    \left( \int_{t_i}^{t_{i + 1}} \int_{A_j} h (t, x) \mathbf{W} (dt, dx)
    \right)^2 \right].
  \end{align*}
Note that
  \begin{align*}
    &\hat{\mathbb{E}} \left[ \sum_{i = 0}^{n - 1} \sum_{j = 1}^m X_{ij}^2
    \left( \int_{t_i}^{t_{i + 1}} \int_{A_j} h (t, x) \mathbf{W} (dt, dx)
    \right)^2 \right] \\
    &\leq \hat{\mathbb{E}} \left[\zeta\right]+\overline{\sigma}^2\hat{\mathbb{E}} \left[ \sum_{i = 0}^{n - 1} \sum_{j = 1}^m X_{ij}^2\int_{t_i}^{t_{i + 1}} \int_{A_j} |h (t, x) |^2 dxdt\right],
  \end{align*}
  where 
   \begin{align*}
    &\hat{\mathbb{E}} \left[\zeta\right]\\
     &= \hat{\mathbb{E}} \left[ \sum_{i = 0}^{n - 1} \sum_{j = 1}^m X_{ij}^2
    \left(\left( \int_{t_i}^{t_{i + 1}} \int_{A_j} h (t, x) \mathbf{W} (dt, dx)
        \right)^2 -\overline{\sigma}^2\int_{t_i}^{t_{i + 1}} \int_{A_j} |h (t, x) |^2 dxdt\right)\right]\\
        &\leq \sum_{i = 0}^{n - 1} \sum_{j = 1}^m\hat{\mathbb{E}} \left[  X_{ij}^2
    \left\{\hat{\mathbb{E}}\left[\left.\left(\int_{t_i}^{t_{i + 1}} \int_{A_j} h (t, x) \mathbf{W} (dt, dx)
            \right)^2\right|\mathcal{F}_{t_i}\right] -\overline{\sigma}^2\int_{t_i}^{t_{i + 1}} \int_{A_j} |h (t, x) |^2 dxdt\right\}\right]\\
        &\leq \sum_{i = 0}^{n - 1} \sum_{j = 1}^m\hat{\mathbb{E}} \left[  X_{ij}^2
    \left(\overline{\sigma}^2\int_{t_i}^{t_{i + 1}} \int_{A_j} |h (t, x)|^2 dxdt -\overline{\sigma}^2\int_{t_i}^{t_{i + 1}} \int_{A_j} |h (t, x) |^2 dxdt\right)\right]\\
        &=0.
  \end{align*}
  Thus,
  \begin{eqnarray*}
    & \hat{\mathbb{E}} \left[ \left| \int_0^T \int_{\mathbb{R}^d} h (t, x) \eta
    (t, x) \mathbf{W} (dt, dx) \right|^2 \right] & \leq \overline{\sigma}^2 
    \hat{\mathbb{E}} \left[ \sum_{i = 0}^{n - 1} \sum_{j = 1}^m X_{ij}^2 
    \int_{t_i}^{t_{i + 1}} \int_{A_j} |h (t, x) |^2 dxdt \right]\\
    &  & \leq \overline{\sigma}^2 \sup_{0 \leq i \leq n - 1} \sup_{1 \leq j \leq
    m}  \hat{\mathbb{E}} [X_{ij}^2]  \int_0^T \int_{\mathbb{R}^d} |h (t, x)
    |^2 dxdt\\
    &  & = \overline{\sigma}^2 \|\eta\|_{\mathbf{S}^2}^2 \|h\|_{L^2}^2 .
  \end{eqnarray*} 
\end{proof}

From the above lemma, the stochastic integral $I_h :
\mathbf{M}^{2, 0} ([0, T] \times \mathbb{R}^d) \rightarrow \mathbf{L}_G^2
(\mathcal{F}_T)$ in Eq. (\ref{stochastic integral2}) can be continuously extended to the mapping from
$\mathbf{S}_G^2 ([0, T] \times \mathbb{R}^d)$ to $\mathbf{L}_G^2
(\mathcal{F}_T)$. Indeed, for any $\eta \in \mathbf{S}_G^2 ([0, T] \times
\mathbb{R}^d)$, there exists a sequence of simple random fields $\eta_n \in
\mathbf{M}^{2, 0} ([0, T] \times \mathbb{R}^d)$ such that $\lim_{n \rightarrow
\infty} \|\eta_n - \eta\|_{\mathbf{S}^2} = 0$. By  inequality (\ref{neq-1}), $\left\{ \int_0^T
\int_{\mathbb{R}^d} h (t, x) \eta_n (t, x) \mathbf{W} (dt, dx)\right\}_{n =
1}^{\infty}$ is a Cauchy sequence in $\mathbf{L}_G^2 (\mathcal{F}_T)$.
Thus we can define
\begin{equation*} 
\int_0^T \int_{\mathbb{R}^d} h (t, x) \eta (t, x) \mathbf{W} (dt, dx) :=
   \mathbb{L}^2 - \lim_{n \rightarrow \infty} \int_0^T \int_{\mathbb{R}^d} h
   (t, x) \eta_n  (t, x) \mathbf{W} (dt, dx) . 
\end{equation*}

{Following the similar procedure as above, for any given set $K\in\mathcal{B}_0(\mathbb{R}^d)$, $h\in L^2([0,T]\times K)$ and $\eta\in\mathbf{S}_G^2([0,T]\times K)$, we can also define the related Bochner integral from $\mathbf{S}_G^2([0,T]\times K)$ to $\mathbf{L}_G^2 (\mathcal{F}_T)$ as follows:
 \begin{align} \int_0^T \int_{K} h (t, x) \eta (t, x)  dxdt&:= \mathbb{L}^2 - \lim_{n \rightarrow \infty} \int_0^T \int_{K} h
   (t, x) \eta_n  (t, x) dxdt \label{b-int1}\\
   &=\mathbb{L}^2 - \lim_{n \rightarrow \infty}\sum_{i = 0}^{m_n - 1} \sum_{j = 1}^{k_n} X_{ij}^{(n)}  \int_{t^{(n)}_i}^{t^{(n)}_{i + 1}}
  \int_{A^{(n)}_j} h(t, x)dxdt,\label{b-int2}
 \end{align}
where $\eta_n(t,x)=\sum_{i = 0}^{m_n - 1}\sum_{j = 1}^{k_n} X_{ij}^{(n)}   \textbf{1}_{[t^{(n)}_i, t^{(n)}_{i + 1})} (t)  \textbf{1}_{A^{(n)}_j} (x) \in
\mathbf{M}^{2, 0} ([0, T] \times K)$ satisfies that $\eta_n$ converges to $\eta$ in $\mathbf{S}_G^2([0,T]\times K)$ as $n\rightarrow\infty$.}

Additionally,  we mention that the Bochner integrals in Eqs. (\ref{bointeg}) and (\ref{bointeg-2}) can also be generalized to the case  that $X_{ij}\in\mathbf{L}_G^p (\mathcal{F}_T)$. This generalization will be utilized in the study of  the stochastic Fubini theorem in the following section. Indeed,
let $\widetilde{\mathbf{M} }^{p, 0} ([0, T] \times K)$ be the
collection of simple random fields with the form:
\begin{equation}
  \eta (s, x ; \omega) = \sum_{i = 0}^{n - 1} \sum_{j = 1}^m X_{ij} (\omega) 
  \textbf{1}_{A_j} (x)  \textbf{1}_{[t_i, t_{i + 1})} (s), \label{simple-func 3}
\end{equation}
where $X_{ij} \in \mathbf{L}_G^p (\mathcal{F}_T)$, $i = 0, \cdots, n -
1$, $j = 1, \cdots, m$, $0 = t_0 < t_1 < \cdots < t_n = T$, and $\{A_j \}_{j =
1}^m \subset \mathcal{B}_0 (\mathbb{R}^d)$ is a mutually disjoint sequence with  $\cup_{j=1}^mA_j\subset K$. Same to definition (\ref{bointeg}), the  Bochner integral for $\eta \in
\widetilde{\mathbf{M} }^{p, 0} ([0, T] \times K)$ can be defined by
\begin{equation*} 
\int_0^T \int_{K} \eta (s, x) dx ds :=
   \sum_{i = 0}^{n - 1} \sum_{j = 1}^m X_{{ij}} (t_{i + 1} - t_i)
   \lambda_{A_j}. 
\end{equation*}
Again, if $K\in\mathcal{B}_0(\mathbb{R}^d)$, this integral  can be continuously extended
to the complete space $\widetilde{\mathbf{M} }^{2}_G ([0, T] \times K)$, where $\widetilde{\mathbf{M} }^{2}_G ([0, T] \times K)$ denotes 
the completion of $\widetilde{\mathbf{M} }^{2, 0} ([0, T] \times
K)$ under the norm $\| \cdot \|_{\widetilde{\mathbf{M}
}^{2}}:= \left( \hat{\mathbb{E}}
   \left[ \int_0^T \int_K |\cdot |^2 dxds \right] \right)^{\frac{1}{2}}$.

\section{Stochastic Fubini Theorem}
In this section, within the sublinear expectation framework, we present a generalization of  the classical  stochastic Fubini theorem in Walsh \cite{Walsh}. This extension  will be applied in our exploration of   weak solutions  for stochastic heat equations.

For any given $T \geq 0$, let $\mathbf{M}^{2, 0} ([0, T]^2 \times \mathbb{R}^2)$ be
the collection of simple random fields in the following form
\begin{equation}
  \label{simple func 2} \eta (t, x, s, y ; \omega) = \sum_{l = 0}^{n_2 - 1} \sum_{k
  = 1}^{m_2} \sum_{i = 0}^{n_1 - 1} \sum_{j = 1}^{m_1} X_{ijlk} (\omega) 
  \textbf{1}_{[t_i, t_{i + 1})} (t)  \textbf{1}_{[s_l, s_{l + 1})} (s) 
  \textbf{1}_{A_j} (x)  \textbf{1}_{B_k} (y),
\end{equation}
where $X_{ijlk} \in \mathbf{L}^2_G (\mathcal{F}_{t_i})$ for $i = 0,
\cdots, n_1 - 1$, $l = 0,
\cdots, n_2 - 1$,  $j = 1, \cdots, m_1$, and $k = 1, \cdots, m_2$, $\{t_0, \cdots, t_{n_1} \}$ and $\{s_0, \cdots, s_{n_2} \}$ are any given
partitions of $[0, T]$, and $A_j , B_k \in
\mathcal{B}_0 (\mathbb{R})$ for $j=1,\cdots,m_1,k=1,\cdots,m_2$ such that
\begin{equation*}
A_{j}\cap A_{j'}=\emptyset, B_{k}\cap B_{k'}=\emptyset \text{ for } 1\leq j,j'\leq m_1,1\leq k,k'\leq m_2 \text { with } j\neq j', k\neq k'.
\end{equation*}

Denote by $\overline{\mathbf{M}}_G^2 ([0, T]^2 \times \mathbb{R}^2)$ the
completion of $\mathbf{M}^{2, 0} ([0, T]^2 \times \mathbb{R}^2)$ under the
norm
\begin{equation*} 
\|\eta (t, x, s, y)\|_{\overline{\mathbf{M}}^2} := \left( \int_0^T
   \int_{\mathbb{R}} \int_0^T {{\int_{\mathbb{R}}} }  \hat{\mathbb{E}} [|\eta (t,
   x, s, y) |^2] dxdtdyds \right)^{1 / 2}. 
\end{equation*}
\begin{theorem}
  \label{Fubini}For $\eta \in \mathbf{M}^{2, 0} ([0, T]^2 \times
  \mathbb{R}^2)$, we have
  \begin{equation}
    \int_0^T \int_{\mathbb{R}} \left[ \int_0^{T}
    \int_{\mathbb{R}} \eta (t, x, s, y) \mathbf{W} (dt, dx) \right] dyds =
    \int_0^T \int_{\mathbb{R}} \left[ \int_0^T \int_{\mathbb{R}} \eta (t, x, s,
    y) dyds \right] \mathbf{W} (dt, dx) . \label{change order}
  \end{equation}
  Furthermore, for $\eta \in \overline{\mathbf{M}}_G^2 ([0, T]^2 \times
  \mathbb{R}^2)$, equality (\ref{change order}) still holds on any compact set
  $K \in \mathcal{B}_0 (\mathbb{R})$, that is,
  \begin{equation}
    \int_0^T \int_K \left[ \int_0^T \int_{\mathbb{R}} \eta (t, x, s, y)
    \mathbf{W} (dt, dx) \right] dyds = \int_0^T \int_{\mathbb{R}} \left[
    \int_0^T \int_K \eta (t, x, s, y) dyds \right] \mathbf{W} (dt, dx) .
    \label{change order2}
  \end{equation}
\end{theorem}

\begin{proof}
  For each $\eta \in \mathbf{M}^{2, 0} ([0, T]^2 \times \mathbb{R}^2)$ with the
  form (\ref{simple func 2}), we have
  \begin{align*}
    & \int_0^T \int_{\mathbb{R}} \eta (t, x, s, y) \mathbf{W} (dt, dx) = \sum_{l
    = 0}^{n_2 - 1} \sum_{k = 1}^{m_2} \left[ \sum_{i = 0}^{n_1 - 1} \sum_{j = 1}^{m_1}
    X_{ijlk}  \mathbf{W} ([t_i, t_{i + 1}) \times A_j) \right] 
    \textbf{1}_{[s_l, s_{l + 1})} (s)  \textbf{1}_{B_k} (y),
  \end{align*}
  and
  \begin{align*}
    & \int_0^T \int_{\mathbb{R}} \eta (t, x, s, y) \textbf{} dyds = \sum_{i =
    0}^{n_1 - 1} \sum_{j = 1}^{m_1} \left[ \sum_{l = 0}^{n_2 - 1} \sum_{k = 1}^{m_2}
    X_{ijlk} (s_{l + 1} - s_l) \lambda_{B_k} \right]  \textbf{1}_{[t_i, t_{i +
    1})} (t)  \textbf{1}_{A_j} (x) \textbf{} .
  \end{align*}
  Obviously, 
  \begin{align*}
  &\left\{\int_0^T \int_{\mathbb{R}} \eta (t, x, s, y) \mathbf{W} (dt, dx)\right\}_{s\in [0,T], y\in \mathbb{R}} \in
  \widetilde{\mathbf{M} }^{2, 0} ([0, T] \times \mathbb{R}), \\ 
  &\left\{\int_0^T
  \int_{\mathbb{R}} \eta (t, x, s, y)  dyds\right\}_{t\in [0,T], x\in \mathbb{R}} \in \mathbf{M}^{2, 0} ([0,
  T]  \times \mathbb{R}),
  \end{align*}
   and the integrals on the both sides of  Eq. (\ref{change order}) are well-defined. Then
  \begin{align*}
      &\int_0^T \int_{\mathbb{R}} \left[ \int_0^T \int_{\mathbb{R}} \eta (t,
    x, s, y) \mathbf{W} (dt, dx) \right] dyds\\
     & = \sum_{l = 0}^{n_2 - 1} \sum_{k = 1}^{m_2} \left[ \sum_{i = 0}^{n_1 - 1}
    \sum_{j = 1}^{m_1} X_{ijlk}  \mathbf{W} ([t_i, t_{i + 1}) \times A_j) \right]
    (s_{l + 1} - s_l) \lambda_{B_k}\\
     & = \int_0^T \int_{\mathbb{R}} \left[ \int_0^T \int_{\mathbb{R}} \eta (t, x, s, y)
    dyds \right] \mathbf{W} (dt, dx),
  \end{align*}
  which implies Eq. (\ref{change order}).
  
  Next for each $\eta \in \overline{\mathbf{M}}_G^2 ([0, T]^2 \times
  \mathbb{R}^2)$,  there exists a sequence
  of simple random fields $\{ \eta_n \}_{n = 1}^{\infty} \subset \mathbf{M}^{2, 0}
  ([0, T]^2 \times \mathbb{R}^2)$ such that $\eta_n \rightarrow \eta$ under the norm $\|\cdot\|_{\overline{\mathbf{M}}^2}$
   as $n\rightarrow\infty$. First we show that the integrals on both sides of Eq. (\ref{change order2}) make sense. For any given $K \in \mathcal{B}_0 (\mathbb{R})$, since
  \begin{align*}
    &\left\| \mathbf{1}_K(y)\int_0^T \int_{\mathbb{R}} \eta_n (t, x, s, y) \mathbf{W} (dt, dx)
    \right\|_{\widetilde{\mathbf{M} }^{2}}^2 \\
    & =  \hat{\mathbb{E}} \left[
    \int_0^T \int_K \left( \int_0^T \int_{\mathbb{R}} \eta_n (t, x, s, y)
    \mathbf{W} (dt, dx) \right)^2 dyds \right]\\
    & \leq  \overline{\sigma}^2 \int_0^T \int_K \hat{\mathbb{E}} \left[ \int_0^T
    \int_{\mathbb{R}} |\eta_n (t, x, s, y) |^2 dxdt \right] dyds\\
    & \leq  \overline{\sigma}^2 \|\eta_n \|_{\overline{\mathbf{M}}^2}^2,
  \end{align*}
  it follows that
  $\left\{ \mathbf{1}_K(y)\int_0^T \int_{\mathbb{R}} \eta_n (t, x, s, y) \mathbf{W} (dt, dx)
  \right\}_{n = 1}^{\infty} \subset \widetilde{\mathbf{M} }^{2, 0} ([0, T]
  \times K)$ is a Cauchy sequence under the norm $\| \cdot
  \|_{\widetilde{\mathbf{M} }^{2}}$ and  converges to  $\mathbf{1}_K(y)\int_0^T
  \int_{\mathbb{R}} \eta (t, x, s, y) \mathbf{W} (dt, dx) \in
  \widetilde{\mathbf{M} }^{2}_G ([0, T] \times K)$. Thus, the
  integral on the left-hand side of Eq. (\ref{change order2}) is well-defined.
  Similarly, by H\"older's inequality,
   \begin{eqnarray*}
       {\left\| \int_0^T \int_K \eta_n (t, x, s, y) d y d s
       \right\|_{\mathbf{M}^2}^2} 
       & \leq &  T \lambda_K \hat{\mathbb{E}} \left[ \int_0^T
       \int_{\mathbb{R}} \int_0^T \int_K | \eta_n (t, x, s, y) |^2 
       dyds d x  dt \right]\\
       & \leq & T \lambda_K \|\eta_n \|_{\overline{\mathbf{M}}^2}^2.
  \end{eqnarray*}
Then  $\left\{ \int_0^T \int_K \eta_n (t, x, s, y) dy d s \right\}_{n = 1}^{\infty}
  \subset \mathbf{M}^{2, 0} ([0, T] \times \mathbb{R})$ is a Cauchy sequence
  under the norm ${\| \cdot \|_{\mathbf{M}^2}}$ and converges to
   $\int_0^T \int_K \eta (t, x, s, y) dyds \in
  \mathbf{M}^2_G ([0, T] \times \mathbb{R})$. Thus, the integral on the right-hand
  side of Eq. (\ref{change order2}) is well-defined.
  
By Lemma \ref{stointprop} and H\"older's inequality, we obtain 
  \begin{align*}
    &\int_0^T \int_K \left( \int_0^T \int_{\mathbb{R}} \eta (t, x, s, y)
    \mathbf{W} (dt, dx) \right) dyds\\
   & =  \mathbb{L}^2 - \lim_{n \rightarrow \infty} \int_0^T \int_K \left(
    \int_0^T \int_{\mathbb{R}} \eta_n (t, x, s, y) \mathbf{W} (dt, dx) \right)
    dyds,
  \end{align*}
  and
  \begin{align*}
   &\int_0^T \int_{\mathbb{R}} \left( \int_0^T \int_K \eta (t, x, s, y)
    dyds \right) \mathbf{W} (dt, dx)\\
    & =  \mathbb{L}^2 - \lim_{n \rightarrow \infty} \int_0^T
    \int_{\mathbb{R}} \left( \int_0^T \int_K \eta_n (t, x, s, y) dyds \right)
    \mathbf{W} (dt, dx) .
  \end{align*}
  As $\eta_n(t,x,s,y)\mathbf{1}_K(y)$ satisfies Eq. (\ref{change order}) for any $n \geq 1$, we have, in $\mathbf{L}_G^2(\mathcal{F}_T)$,
  \begin{equation*} 
  \int_0^T \int_K \left[ \int_0^T \int_{\mathbb{R}} \eta (t, x, s, y)
     \mathbf{W} (dt, dx) \right] dyds = \int_0^T \int_{\mathbb{R}} \left[
     \int_0^T \int_K \eta(t, x, s, y) dyds \right] \mathbf{W} (dt, dx) . 
  \end{equation*}
\end{proof}

Building upon the above theorem, we now present a specific instance of the stochastic Fubini theorem. This result will be employed in proving the existence of  weak solutions for  $G$-stochastic heat equations.

\begin{corollary}\label{Fubi-2}
Suppose $h\in L^2([0,T]\times \mathbb{R})$, $\phi\in C_c^\infty([0,T]\times\mathbb{R})$ and $\eta\in \mathbf{S}_G^2([0,T]\times\mathbb{R})$, where $C_c^{\infty} ([0,T]\times\mathbb{R})$ denotes the space of infinitely
differentiable functions with compact support in $[0,T]\times\mathbb{R}$. For each $0<s,t\leq T$, $x,y\in\mathbb{R}$, set
\begin{align*}
&g(s,y):=\int_s^T\int_\mathbb{R}h(t-s,x-y)\phi(t,x)dxdt,\\
&\zeta(t,x):=\int_0^t\int_\mathbb{R}h(t-s,x-y)\eta(s,y)\mathbf{W}(ds,dy).
\end{align*}
If $\zeta\in \mathbf{S}_G^2([0,T]\times \mathbb{R})$ and
\begin{align}\label{integ-h}
&\int_0^T\int_{\mathbb{R}}\int_s^T\int_\mathbb{R}h^2(t-s,x-y)\phi^2(t,x)dxdtdyds<\infty,
\end{align}
then, in $\mathbf{L}_G^2(\mathcal{F}_T)$,
\begin{align}\label{Fueq-2}
&\int_0^T\int_\mathbb{R}g(s,y)\eta(s,y)\mathbf{W}(ds,dy)=\int_0^T\int_\mathbb{R}\zeta(t,x)\phi(t,x)dxdt.
\end{align}
\end{corollary}
\begin{proof}
By assumption (\ref{integ-h}), it is easy to check that $g\in L^2([0,T]\times \mathbb{R})$. Then the integrals on  both sides of Eq. (\ref{Fueq-2}) make sense. 

For each $\eta\in \mathbf{S}_G^2([0,T]\times\mathbb{R})$, there exists a sequence of simple random fields $\{\eta_n\}_{n=1}^\infty\subset \mathbf{M}^{2,0}([0,T]\times\mathbb{R})$ such that $\|\eta_n-\eta\|_{\mathbf{S}^2}\rightarrow 0$ as $n\rightarrow\infty$. Assume that 
\begin{equation*}
  \eta_n (s, y) = \sum_{i = 0}^{k_n - 1} \sum_{j = 1}^{m_n} X^{(n)}_{ij}
  \textbf{1}_{[t^n_i, t^n_{i + 1})} (s) \textbf{1}_{A^n_j} (y) , 
\end{equation*}
where $X^{(n)}_{ij} \in \mathbf{L}_G^2 (\mathcal{F}_{t^n_i})$, $i = 0, \cdots, k_n -
1$, $j = 1, \cdots, m_n$, $0 = t^n_0 < t^n_1 < \cdots < t^n_{k_n} = T$, and $\{A^n_j \}_{j =
1}^{m_n} \subset \mathcal{B}_0 (\mathbb{R}^d)$ is a mutually disjoint sequence. 
Then for $n\geq1$, 
\begin{align*}
&\int_0^T\int_\mathbb{R}g(s,y)\eta_n(s,y)\mathbf{W}(ds,dy)\\
&=\sum_{i = 0}^{k_n - 1} \sum_{j = 1}^{m_n}X^{(n)}_{ij}\int_{t^n_i}^{t^n_{i+1}}\int_{A_j^n}g(s,y)\mathbf{W}(ds,dy)\\
&=\sum_{i = 0}^{k_n - 1} \sum_{j = 1}^{m_n}X^{(n)}_{ij}\int_{t^n_i}^{t^n_{i+1}}\int_{A_j^n}\left(\int_0^T\int_\mathbb{R}h(t-s,x-y)\phi(t,x)\mathbf{1}_{[s,T]}(t)dxdt\right)\mathbf{W}(ds,dy),
\end{align*}
and
\begin{align*}
&\int_0^T\int_\mathbb{R}\zeta_n(t,x)\phi(t,x)dxdt\\
&=\int_0^T\int_\mathbb{R}\left( \sum_{i = 0}^{k_n - 1} \sum_{j = 1}^{m_n} X^{(n)}_{ij}\int_{t^n_i}^{t^n_{i+1}}\int_{A_j^n}h(t-s,x-y)\mathbf{1}_{[0,t]}(s)\mathbf{W}(ds,dy)\right)\phi(t,x)dxdt\\
&= \sum_{i = 0}^{k_n - 1} \sum_{j = 1}^{m_n} X^{(n)}_{ij}\int_0^T\int_\mathbb{R}\left(\int_{t^n_i}^{t^n_{i+1}}\int_{A_j^n}h(t-s,x-y)\phi(t,x)\mathbf{1}_{[0,t]}(s)\mathbf{W}(ds,dy)\right)dxdt,
\end{align*}
where $\zeta_n(t,x):=\int_0^t\int_\mathbb{R}h(t-s,x-y)\eta_n(s,y)\mathbf{W}(ds,dy)$.
By  assumption (\ref{integ-h}), {we know that $h(t-s,x-y)\phi(t,x)\mathbf{1}_{[0,t]}(s)\in L^2([0,T]^2\times\mathbb{R}^2)\subset\overline{\mathbf{M}}_G^2([0,T]^2\times\mathbb{R}^2)$.} Applying Theorem \ref{Fubini},  we deduce that
\begin{equation}\label{equali-1}
\int_0^T\int_\mathbb{R}g(s,y)\eta_n(s,y)\mathbf{W}(ds,dy)=\int_0^T\int_\mathbb{R}\zeta_n(t,x)\phi(t,x)dxdt, \text{ for } n\geq1.
\end{equation}

Now we show  that the integrals in Eq. (\ref{equali-1}) converge in $\mathbf{L}_G^2(\mathcal{F}_T)$. As
\begin{align*}
\hat{\mathbb{E}}\left[\left|\zeta_n(t,x)-\zeta(t,x)\right|^2\right]
&= \hat{\mathbb{E}}\left[\left|\int_0^t\int_\mathbb{R}h(t-s,x-y)(\eta_n(s,y)-\eta(s,y))\mathbf{W}(ds,dy)\right|^2\right]\\
&\leq \overline{\sigma}^2\int_0^t\int_\mathbb{R}|h(t-s,x-y)|^2\hat{\mathbb{E}}\left[|\eta_n(s,y)-\eta(s,y)|^2\right]dyds\\
&\leq \overline{\sigma}^2\left\|\eta_n-\eta\right\|^2_{\mathbf{S}^2}\int_0^t\int_\mathbb{R}|h(t-s,x-y)|^2dyds,
\end{align*}
we have
\begin{align*}
&\hat{\mathbb{E}}\left[\left|\int_0^T\int_\mathbb{R}(\zeta_n(t,x)-\zeta(t,x))\phi(t,x)dxdt\right|^2\right]\\
&\leq C_{T,\phi}\int_0^T\int_\mathbb{R}\hat{\mathbb{E}}\left[\left|\zeta_n(t,x)-\zeta(t,x)\right|^2\right]\left|\phi(t,x)\right|^2dxdt\\
&\leq C_{T,\phi}\overline{\sigma}^2\left\|\eta_n-\eta\right\|^2_{\mathbf{S}^2}\int_0^T\int_\mathbb{R}\int_0^t\int_\mathbb{R}|h(t-s,x-y)|^2\left|\phi(t,x)\right|^2dydsdxdt\\
&\leq C_{T,h,\phi}\overline{\sigma}^2 \left\|\eta_n-\eta\right\|^2_{\mathbf{S}^2},
\end{align*}
where $C_{T,\phi}$ and $C_{T,h,\phi}$ are constants.
For the left-hand side of Eq. (\ref{equali-1}),
\begin{align*}
&\hat{\mathbb{E}}\left[\left|\int_0^T\int_\mathbb{R}g(s,y)(\eta_n(s,y)-\eta(s,y))\mathbf{W}(ds,dy)\right|^2\right]\\
&\leq \overline{\sigma}^2\hat{\mathbb{E}}\left[\int_0^T\int_\mathbb{R}\left|g(s,y)\right|^2\left|\eta_n(s,y)-\eta(s,y)\right|^2dyds\right]\\
&\leq \overline{\sigma}^2 \int_0^T\int_\mathbb{R}\left|g(s,y)\right|^2\hat{\mathbb{E}}\left[\left|\eta_n(s,y)-\eta(s,y)\right|^2\right]dyds\\
&\leq \overline{\sigma}^2 \|g\|^2_{L^2}\left\|\eta_n-\eta\right\|^2_{\mathbf{S}^2}.
\end{align*}
Then taking $n\rightarrow\infty$ in Eq. (\ref{equali-1}), both sides of Eq. (\ref{equali-1}) converge in $\mathbf{L}_G^2(\mathcal{F}_T)$ and 
\begin{align*}
&\int_0^T\int_\mathbb{R}g(s,y)\eta(s,y)\mathbf{W}(ds,dy)=\int_0^T\int_\mathbb{R}\zeta(t,x)\phi(t,x)dxdt.
\end{align*}
\end{proof}

\section{Stochastic heat equations driven by space-time $G$-white noises}

In this section, we focus on the linear and nonlinear stochastic heat equations driven by the 
space-time $G$-white noise.

For each fixed $T > 0$,   consider the
space-time $G$-white noise $\mathbf{W} = \{\mathbf{W}
([s, t) \times A) : 0\leq s\leq t \leq T, A \in \mathcal{B}_0(\mathbb{R}) \}$. Without causing confusion, for simplicity of notation,  denote by 
\begin{equation*}
\mathbf{W} (t, x):=\mathbf{W}
([0, t) \times [0\wedge x, 0\vee x]), \text{ for }t \in [0, T], x \in \mathbb{R}.
\end{equation*}
By the definition and properties of the space-time $G$-white noise, we know that $\mathbf{W} (t, x)$ is
nowhere-differentiable with respect to $t$ and $x$ in the ordinary sense. Nevertheless, we may define its derivative  in the sense of  generalized functions.  For each test function $\phi \in C_c^{\infty}(\mathbb{R}^2)$, as  $(\mathbf{W}(t,x))_{t\in[0,T],x\in K}\in\mathbf{S}_G^2([0,T]\times K)$ for any compact set $K$, the integral 
\begin{equation}
  I_{\mathbf{W}}(\phi):=\int_0^T \int_{\mathbb{R}} \mathbf{W} (t, x) \frac{\partial^2 \phi (t,
  x)}{\partial t \partial x} dxdt \label{W-varphi}
\end{equation}
makes sense and $I_{\mathbf{W}}$ is a linear mapping from $C_c^{\infty}(\mathbb{R}^2)$ to $\mathbf{L}_G^2(\mathcal{F}_T)$.
Define the generalized mixed derivative $\frac{\partial^2 \mathbf{W} (t,
   x)}{\partial t \partial x}$ through the test function  $\phi \in C_c^{\infty}(\mathbb{R}^2)$ as follows 
\begin{equation}
  \int_0^T \int_{\mathbb{R}} \frac{\partial^2 \mathbf{W} (t,
   x)}{\partial t \partial x} \phi (t, x) dxdt :=
  \int_0^T \int_{\mathbb{R}} \mathbf{W} (t, x) \frac{\partial^2 \phi (t,
  x)}{\partial t \partial x} dxdt. \label{deri def}
\end{equation}
Then $\frac{\partial^2 \mathbf{W} (t,
   x)}{\partial t \partial x}$  is called the derivative of ${\mathbf{W}} (t, x)$, denoted by $\dot{\mathbf{W}} (t, x)$, $t \in [0, T], x \in \mathbb{R}$.

By Proposition 2.8 in Chapter 2 in Khoshnevisan {\cite{Kh14}}, it is easy to get the following equality of stochastic integrals.

\begin{proposition}
  \label{inte equal prop}Let $\{ \dot{\mathbf{W}} (t, x) : t \in [0, T], x \in
  \mathbb{R} \}$ be the derivative of the space-time $G$-white noise $\{
  \mathbf{W} (t, x) : t \in [0, T], x \in \mathbb{R} \}$ in the sense of
  (\ref{deri def}). Then for each $\phi \in C_c^{\infty} (\mathbb{R}^2)$, we
  have
  \begin{equation}
    \int_0^T \int_{\mathbb{R}} \dot{\mathbf{W}} (t, x) \phi (t, x) dxdt =
    \int_0^T \int_{\mathbb{R}} \phi (t, x) \mathbf{W} (dt, dx) . \label{integ
    equal}
  \end{equation}
\end{proposition}

\begin{proof}
  By Eq. (\ref{deri def}), we only need to verify that for any $\phi \in
  C_c^{\infty} (\mathbb{R}^2)$, in $\mathbf{L}_G^2(\mathcal{F}_T)$, 
  \begin{equation}\label{inte-equal} 
  \int_0^T \int_{\mathbb{R}} \mathbf{W} (t, x) \frac{\partial^2 \phi (t,
     x)}{\partial t \partial x} dxdt = \int_0^T \int_{\mathbb{R}} \phi (t, x)
     \mathbf{W}(dt, dx) . 
\end{equation}
  Without loss of generality, assume that $T \geq 1$ and  $\phi$ is supported in $(0, 1)
  \times (0, 1)$.
   For each $n \geq1$, set
  \begin{equation} 
  \phi_n (t, x) := \sum_{j = 0}^{n - 1} \sum_{i = 0}^{n - 1} \phi \left(
     \frac{i}{n}, \frac{j}{n} \right)  \textbf{1}_{[\frac{i}{n},
     \frac{i + 1}{n})}(t)\textbf{1}_{ [\frac{j}{n} ,\frac{j + 1}{n} )}(x), \text{ for } t,x\in [0,1]. 
  \end{equation}
Then $\phi_n\in\mathbf{M}^{2,0}([0,1]\times[0,1])$,  and 
  \begin{align*}
    &\int_0^T \int_{\mathbb{R}} \phi_n (t, x) \mathbf{W} (dt, dx)\\
     & = \sum_{j =
    0}^{n - 1} \sum_{i = 0}^{n - 1} \phi \left( \frac{i}{n}, \frac{j}{n}
    \right) \mathbf{W} \left( \left[ \frac{i}{n}, \frac{i + 1}{n} \right)
    \times \left[ \frac{j}{n}, \frac{j + 1}{n} \right) \right)\nonumber\\
    & = \sum_{j = 0}^{n - 1} \sum_{i = 0}^{n - 1} \phi \left( \frac{i}{n},
    \frac{j}{n} \right)  \left[ \mathbf{W} \left( \frac{i + 1}{n}, \frac{j +
    1}{n} \right) - \mathbf{W} \left( \frac{i}{n}, \frac{j + 1}{n} \right) -
    \mathbf{W} \left( \frac{i + 1}{n}, \frac{j}{n} \right) + \mathbf{W} \left(
    \frac{i}{n}, \frac{j}{n} \right) \right]\nonumber\\
    & = \sum_{j = 1}^n \sum_{i = 1}^n \mathbf{W} \left( \frac{i}{n},
    \frac{j}{n} \right)  \left[ \phi \left( \frac{i}{n}, \frac{j}{n} \right) -
    \phi \left( \frac{i - 1}{n}, \frac{j}{n} \right) - \phi \left(
    \frac{i}{n}, \frac{j - 1}{n} \right) + \phi \left( \frac{i - 1}{n},
    \frac{j - 1}{n} \right) \right].
  \end{align*}
By transforming the expression of $\phi$ inside the bracket into integral form, we have
  \begin{align}
    \int_0^T \int_{\mathbb{R}} \phi_n (t, x) \mathbf{W} (dt, dx)
    & = \sum_{j = 1}^{n } \sum_{i = 1}^{n } \mathbf{W} \left(\frac{i}{n}, \frac{j}{n} \right)\int_{\frac{i-1}{n}}^{\frac{i}{n}} \int_{\frac{j-1}{n}}^{\frac{j}{n}}  \frac{\partial^2 \phi (t, x)}{\partial t
    \partial x} dxdt,\ \forall n\geq 1.\label{eqinte}
  \end{align}

Since $\phi_n \rightarrow \phi $ in $L^2 (\mathbb{R}^2)$ as $n
  \rightarrow \infty$, by Lemma \ref{stointprop}, we deduce that
  \begin{equation}\label{L2Wphi-1}
  \int_0^T \int_{\mathbb{R}} \phi (t, x) \mathbf{W}
  (dt, dx) =\mathbb{L}^2-\lim_{n\rightarrow\infty}\int_0^T \int_{\mathbb{R}} \phi_n (t, x) \mathbf{W} (dt,dx).
  \end{equation}
Thus, to obtain the equality (\ref{inte-equal}), it suffices to prove 
\begin{equation}\label{L2Wphi}
\int_0^T \int_{\mathbb{R}} \mathbf{W} (t, x) \frac{\partial^2 \phi (t,x)}{\partial t \partial x} dxdt=\mathbb{L}^2-\lim_{n\rightarrow\infty}\sum_{j = 1}^{n } \sum_{i = 1}^{n } \mathbf{W} \left(\frac{i}{n}, \frac{j}{n} \right)\int_{\frac{i-1}{n}}^{\frac{i}{n}} \int_{\frac{j-1}{n}}^{\frac{j}{n}}  \frac{\partial^2 \phi (t, x)}{\partial t \partial x} dxdt.
  \end{equation} 
Indeed, for any $n\geq1$, 
  \begin{align*}
     & \hat{\mathbb{E}} \left[ \left|
  \int_0^T \int_{\mathbb{R}} \mathbf{W} (t, x) \frac{\partial^2 \phi (t,x)}{\partial t \partial x} dxdt-\sum_{j = 1}^{n } \sum_{i = 1}^{n } \mathbf{W} \left(\frac{i}{n}, \frac{j}{n} \right)\int_{\frac{i-1}{n}}^{\frac{i}{n}} \int_{\frac{j-1}{n}}^{\frac{j}{n}}  \frac{\partial^2 \phi (t, x)}{\partial t
    \partial x} dxdt\right|^2 \right]\\
     & =\hat{\mathbb{E}} \left[ \left|\sum_{j = 1}^{n } \sum_{i = 1}^{n }\int_{\frac{i-1}{n}}^{\frac{i}{n}} \int_{\frac{j-1}{n}}^{\frac{j}{n}}\left(\mathbf{W} (t, x)-\mathbf{W} \left(\frac{i}{n}, \frac{j}{n} \right)\right)  \frac{\partial^2 \phi (t,x)}{\partial t \partial x} dxdt\right|^2 \right]\\
      & \leq\hat{\mathbb{E}} \left[ \left(\sum_{j = 1}^{n } \sum_{i = 1}^{n } \int_{\frac{i-1}{n}}^{\frac{i}{n}}\int_{\frac{j-1}{n}}^{\frac{j}{n}}\left|\mathbf{W} (t, x)-\mathbf{W} \left(\frac{i}{n}, \frac{j}{n} \right)\right|^2  dxdt\right) \left(\int_{0}^{1}\int_{0}^{1}  \left|\frac{\partial^2 \phi (t,x)}{\partial t \partial x}\right| ^2dxdt\right) \right]\\
       & \leq C_\phi\sum_{j = 1}^{n } \sum_{i = 1}^{n } \int_{\frac{i-1}{n}}^{\frac{i}{n}}\int_{\frac{j-1}{n}}^{\frac{j}{n}}\hat{\mathbb{E}} \left[ \left|\mathbf{W} (t, x)-\mathbf{W} \left(\frac{i}{n}, \frac{j}{n} \right)\right|^2\right]  dxdt \\
    & \leq C_\phi\overline{\sigma}^2  \sum_{j = 1}^{n } \sum_{i = 1}^{n } \int_{\frac{i-1}{n}}^{\frac{i}{n}}\int_{\frac{j-1}{n}}^{\frac{j}{n}}
    \left(\frac{i j}{n^2}-tx \right) dxdt\\
     & = \frac{ C_\phi\overline{\sigma}^2(2n+1)}{4n^2} ,
  \end{align*}
where $C_\phi$ is a constant depending on $\phi$. Letting $n \rightarrow \infty$ in the above inequality, we get Eq. (\ref{L2Wphi}) and complete the proof.
\end{proof}

\subsection{Linear $G$-stochastic heat equations}

In this subsection, we focus on the linear stochastic heat equation driven by the additive space-time $G$-white noise as follows:
\begin{equation}
  \left\{\begin{array}{ll}
    & \frac{\partial}{\partial t} u (t, x) - \frac{\partial^2}{\partial x^2}
    u (t, x) = \dot{\mathbf{W}} (t, x),\  t > 0, x \in \mathbb{R},\\
    & u (0, x) = u_0 (x), x \in \mathbb{R},
  \end{array}\right. \label{linear equation}
\end{equation}
where the initial function $u_0 : \mathbb{R} \rightarrow \mathbb{R}$ is  nonrandom, bounded and Borel measurable.

Applying Proposition \ref{inte equal prop}, we can define the weak solution of  (\ref{linear equation}).

\begin{definition}\label{weak-sol-def1}
  A spatio-temporal random field $\{u (t, x) : (t, x) \in [0, \infty) \times
  \mathbb{R} \}$ is said to be a  {weak solution} of the linear $G$-stochastic heat equation
  (\ref{linear equation}) if it satisfies the following conditions:
  \begin{itemize}
    \item For each $T > 0$ and compact set $K \in \mathcal{B}_0 (\mathbb{R})$,
    {$(u (t, x))_{t\in[0,T], x \in K} \in \mathbf{S}_G^2 ([0, T]
        \times K)$}; 
    
    \item For each $T > 0$ and $\phi  \in C_c^{\infty}
    (\mathbb{R}^2)$, 
    \begin{eqnarray}
      &  & \int_{\mathbb{R}} u (T, x) \phi (T, x) dx - \int_{\mathbb{R}} u_0
      (x) \phi (0, x) dx \label{linear weak def}\\
      &  & = \int_0^T \int_{\mathbb{R}} u (t, x)  \left( \frac{\partial \phi
      (t, x)}{\partial t} + \frac{\partial^2 \phi (t, x)}{\partial x^2}
      \right) dxdt + \int_0^T \int_{\mathbb{R}} \phi (t, x) \mathbf{W} (dt,
      dx) .  \nonumber
    \end{eqnarray}
  \end{itemize}
\end{definition}

Now we present the explicit expression for the weak solution of  (\ref{linear equation}).

\begin{theorem}
  \label{linear eq thm}For each given bounded and Borel measurable function $u_0 (x)$, the linear
  $G$-stochastic heat equation (\ref{linear equation})   has a weak solution defined by
\begin{equation}
  u (t, x) = \int_{\mathbb{R}} u_0 (y) p (t, x - y) dy + \int_0^t
  \int_{\mathbb{R}} p (t - s, x - y) \mathbf{W} (ds, dy),\ \text{for } t > 0, x \in
  \mathbb{R}, \label{linear solution}
\end{equation}
and $u(0,x)=u_0(x)$ for $x\in\mathbb{R}$. Here $p (t, x) = \frac{1}{\sqrt{4 \pi t}} e^{- \frac{x^2}{4 t}}, t > 0, x \in
\mathbb{R}$. 

Furthermore, if $u,v$ are two weak solutions of  $G$-stochastic heat equation  (\ref{linear equation}) with the same initial function, then  $u(t,x)=v(t,x)$ for almost everywhere $(t,x)\in[0, \infty) \times
  \mathbb{R}$.
\end{theorem}
\begin{remark}
The solution (\ref{linear solution}) is also called the mild solution of $G$-stochastic heat equation (\ref{linear equation}).
\end{remark}
\begin{remark}
  Note that $p (t, x) \in L^2 ([0, T] \times \mathbb{R}),\ \forall T > 0.$ By
  the definition of stochastic integrals in Section 3, the right-hand side of Eq. (\ref{linear
  solution}) is well-defined. In particular, $p (t, x)$ is called the heat kernel (or fundamental solution) and is the solution of the
  following heat equation:
  \begin{equation}
    \left\{\begin{array}{ll}
      & \frac{\partial}{\partial t} p (t, x) - \frac{\partial^2}{\partial
      x^2} p (t, x) = 0, \ t > 0, x \in \mathbb{R},\\
      & \lim_{t \downarrow 0} p (t, x) = \delta (x),\ x \in \mathbb{R} .
    \end{array}\right. \label{funde sol}
  \end{equation}
\end{remark}

To prove Theorem \ref{linear eq thm}, we first provide the moment estimates for the stochastic integral in the definition of weak solution (\ref{linear solution}).

\begin{lemma}\label{pW-est}
Denote
\begin{equation*}
Z(t,x):=\int_0^t\int_{\mathbb{R}} p (t - s, x - y) \mathbf{W} (ds, dy), \ t>0,x\in\mathbb{R},
\end{equation*} 
and $Z(0,x):=0$, $x\in\mathbb{R}$.
Then there exists a constant $C$ depending on $\overline{\sigma}^2$ such that for each $\delta\geq0$, 
\begin{itemize}
    \item[(i)] $\hat{\mathbb{E}} [|Z(t, x + \delta) - Z(t, x) |^2] \leq
    C\delta ;$
    \item[(ii)] $\hat{\mathbb{E}} [|Z(t + \delta, x) - Z(t, x) |^2] \leq C
    \sqrt{\delta} .$
  \end{itemize}
\end{lemma}
\begin{proof}
(i) For each $\delta\geq0$, by Eq. (3.16) in Khoshnevisan \cite{Kh14},  
\begin{align}\label{ineq-p-x}
\int_0^t\int_{\mathbb{R}} \left|p (t - s, x+\delta - y)-p (t - s, x - y)\right|^2 dyds&\leq\frac{\delta}{2}.
\end{align}
Then 
\begin{align*}
&\hat{\mathbb{E}}\left[\left|Z(t, x + \delta) - Z(t, x) \right|^2\right]\\
&=\hat{\mathbb{E}}\left[\left|\int_0^t\int_{\mathbb{R}} (p (t - s, x+\delta - y)-p (t - s, x - y)) \mathbf{W} (ds, dy)\right|^2\right]\\
&\leq \overline{\sigma}^2\int_0^t\int_{\mathbb{R}} \left|p (t - s, x+\delta - y)-p (t - s, x - y)\right|^2 dyds\\
&\leq\frac{\overline{\sigma}^2}{2}\delta.
\end{align*}
(ii) For $\delta\geq0$, since 
\begin{align}\label{ineq-p-t1}
\int_0^t\int_{\mathbb{R}} \left|p (t +\delta- s, x- y)-p (t - s, x - y)\right|^2 dyds\leq(\sqrt{2}-1)\sqrt{\frac{\delta}{2\pi}},
\end{align}
and
\begin{align}\label{ineq-p-t2}
\int_t^{t+\delta}\int_{\mathbb{R}} \left|p (t +\delta- s, x- y)\right|^2 dyds=\sqrt{\frac{\delta}{2\pi}},
\end{align}
we have
\begin{align*}
&\hat{\mathbb{E}}\left[\left|Z(t+\delta, x) - Z(t, x) \right|^2\right]\\
&\leq\hat{\mathbb{E}}\left[\left|\int_0^t\int_{\mathbb{R}} (p (t+\delta - s, x- y)-p (t - s, x - y)) \mathbf{W} (ds, dy)\right|^2\right]\\
&\ \ \ +\hat{\mathbb{E}}\left[\left|\int_t^{t+\delta}\int_{\mathbb{R}} p (t+\delta - s, x- y) \mathbf{W} (ds, dy)\right|^2\right]\\
&\leq \overline{\sigma}^2\int_0^t\int_{\mathbb{R}} \left|p (t +\delta- s, x- y)-p (t - s, x - y)\right|^2 dyds\\
&\ \ \ + \overline{\sigma}^2\int_t^{t+\delta}\int_{\mathbb{R}} \left|p (t +\delta- s, x- y)\right|^2 dyds\\
&\leq\overline{\sigma}^2\sqrt{\frac{\delta}{\pi}}.
\end{align*}
\end{proof}

\begin{proof}[Proof of Theorem \ref{linear eq thm}]
  First, we prove that the random field $(u (t, x))_{t \geq0, x \in \mathbb{R}}$
  defined by Eq. (\ref{linear solution}) satisfies (i) in Definition \ref{weak-sol-def1}. As $u_0$ is a bounded function, it is easy to check that $\int_{\mathbb{R}} u_0 (y) p (t, x - y) dy$ is bounded for  $t\geq0, x\in\mathbb{R}$ and $\left(\int_{\mathbb{R}} u_0 (y) p (t, x - y) dy\right)_{t\in[0,T],x\in K}\in \mathbf{S}_G^2 ([0, T]\times K)$ for any compact set $K\in\mathcal{B}_0(\mathbb{R})$. 

  For any given $T>0$ and compact set $K\in\mathcal{B}_0(\mathbb{R})$, it remains to show 
  \begin{equation*}
  \left\{\int_0^t\int_{\mathbb{R}} p (t - s, x - y) \mathbf{W} (ds, dy): t\in(0,T],x\in K\right\}\in\mathbf{S}_G^2 ([0, T]\times K).
  \end{equation*}
   Without loss of generality, assume that $K\subset[-L,L]$ for some $L>0$. Then it suffices to verify that $\left\{\int_0^t\int_{\mathbb{R}} p (t - s, x - y) \mathbf{W} (ds, dy): t\in(0,T],x\in [-L,L]\right\}\in \mathbf{S}_G^2 ([0, T]\times [-L, L])$.

  For $t\in(0,T]$, $x\in [-L,L]$, $n\geq1$,  denote
  \begin{align*}
  &Z(t,x)=\int_0^t\int_{\mathbb{R}} p (t - s, x - y) \mathbf{W} (ds, dy), \\
  &Z_n(t,x):=\sum_{i=0}^{n-1}\sum_{j=-n}^{n-1}\int_{0}^{\frac{iT}{n}}\int_{\mathbb{R}} p \left(\frac{iT}{n} - s, \frac{jL}{n} - y\right) \mathbf{W} (ds, dy)\mathbf{1}_{[\frac{iT}{n},\frac{(i+1)T}{n})}(t)\mathbf{1}_{[\frac{jL}{n},\frac{(j+1)L}{n})}(x),
  \end{align*}
and $Z(0,x)=Z_n(0,x)=0$, for $x\in[-L,L]$, $n\geq1$.  Then {$\{Z_n\}_{n=1}^\infty\subset \mathbf{M}^{2,0}([0,T]\times[-L,L]) $}. By Lemma \ref{pW-est}, there exists a constant $C$ depending on $T$, $L$ and $\overline{\sigma}^2$ such that
  \begin{align*}
  &\sup_{t\in[0,T]}\sup_{x\in[-L,L]}\hat{\mathbb{E}}[|Z(t,x)-Z_n(t,x)|^2]\leq \frac{C}{\sqrt{n}}.
  \end{align*}
Hence,
\begin{equation*}
\lim_{n\rightarrow\infty}\|Z-Z_n\|_{\mathbf{S}^2}^2=0,
\end{equation*}
and $(Z(t, x))_{t \in[0,T], x \in [-L,L]} \in{\mathbf{S}_G^2 ([0, T] \times [-L,L])}$, which implies   $(u(t, x))_{t \in[0,T], x \in K} \in
    {\mathbf{S}_G^2 ([0, T] \times K)}$ for any $T>0$ and compact set $K\in\mathcal{B}_0(\mathbb{R})$.
  
 Next we show that $u (t, x)$ defined by (\ref{linear solution}) satisfies Eq. (\ref{linear weak def}). Substituting Eq. (\ref{linear solution}) into the first term on the right-hand side
  of Eq. (\ref{linear weak def}), we have
  \begin{eqnarray}
    &  & \int_0^T \int_{\mathbb{R}} u (t, x)  \left( \frac{\partial \phi (t,
    x)}{\partial t} + \frac{\partial^2 \phi (t, x)}{\partial x^2} \right)
    dxdt  \nonumber\\ 
    &  & = \int_0^T \int_{\mathbb{R}} \left( \int_{\mathbb{R}} u_0 (y) p (t,
    x - y) dy \right)  \left( \frac{\partial \phi (t, x)}{\partial t} +
    \frac{\partial^2 \phi (t, x)}{\partial x^2} \right) dxdt  \nonumber\\
    &  & \,  \hspace{1.2em} + \int_0^T \int_{\mathbb{R}} \left( \int_0^t
    \int_{\mathbb{R}} p (t - s, x - y) \mathbf{W} (ds, dy) \right)  \left(
    \frac{\partial \phi (t, x)}{\partial t} + \frac{\partial^2 \phi (t,
    x)}{\partial x^2} \right) dxdt  \nonumber\\
    &  & \triangleq \mathcal{I}_1 +\mathcal{I}_2 .\label{I12}
  \end{eqnarray}
  By the Fubini theorem and integration by parts, we get
  \begin{eqnarray*}
    \mathcal{I}_1 & = & \int_{\mathbb{R}} u_0 (y) \left[ \int_0^T
    \int_{\mathbb{R}} p (t, x - y) \left( \frac{\partial \phi (t, x)}{\partial
    t} + \frac{\partial^2 \phi (t, x)}{\partial x^2} \right) dxdt \right] dy\\
    & = & \int_{\mathbb{R}} \int_{\mathbb{R}} u_0 (y) p (T, x - y) \phi (T,
    x) dxdy - \int_{\mathbb{R}} u_0 (y) \phi (0, y) dy.
  \end{eqnarray*}
  Applying Corollary \ref{Fubi-2}, 
  \begin{eqnarray*}{}
    \mathcal{I}_2 & = & \int_0^T \int_{\mathbb{R}} \left[ \int_s^T
    \int_{\mathbb{R}} p (t - s, x - y) \left( \frac{\partial \phi (t,
    x)}{\partial t} + \frac{\partial^2 \phi (t, x)}{\partial x^2} \right) dxdt
    \right] \mathbf{W} (ds, dy)\\
    & = & \int_0^T \int_{\mathbb{R}} \left[ \int_{\mathbb{R}} p (T - s, x -
    y) \phi (T, x) dx \right] \mathbf{W} (ds, dy) - \int_0^T \int_{\mathbb{R}}
    \phi (s, y) \mathbf{W} (ds, dy) .
  \end{eqnarray*}
  Substituting $\mathcal{I}_1$ and $\mathcal{I}_2$ into Eq.
  (\ref{I12}), it follows that  
  \begin{align}
    & \int_0^T \int_{\mathbb{R}} u (t, x)  \left( \frac{\partial \phi (t,
    x)}{\partial t} + \frac{\partial^2 \phi (t, x)}{\partial x^2} \right) dxdt
    + \int_0^T \int_{\mathbb{R}} \phi (t, x) \mathbf{W} (dt, dx) \nonumber\\
    &  = \int_{\mathbb{R}} \int_{\mathbb{R}} u_0 (y) p (T, x - y) \phi (T,
    x) dxdy - \int_{\mathbb{R}} u_0 (y) \phi (0, y) dy \label{I-2}\\
    &\ \ \  + \int_0^T \int_{\mathbb{R}} \left[
    \int_{\mathbb{R}} p (T - s, x - y) \phi (T, x) dx \right] \mathbf{W} (ds,
    dy) .\nonumber 
  \end{align}
  At $t=T$,  multiplying $\phi(T,x)$ and integrating  both sides of Eq. (\ref{linear solution}) with respect to $x$, we obtain
  \begin{align}
      \int_{\mathbb{R}} u (T, x) \phi (T, x) dx  = &\int_{\mathbb{R}} \int_{\mathbb{R}} u_0 (y) p (T, x - y) \phi (T,x) dxdy \label{neq3}\\
      &+ \int_0^T \int_{\mathbb{R}} \left[ \int_{\mathbb{R}} p (T - s, x
    - y) \phi (T, x) dx \right] \mathbf{W} (ds, dy).  \nonumber
  \end{align}
  Substituting Eq. (\ref{neq3}) into the right-hand side of Eq. (\ref{I-2}), we get   Eq. (\ref{linear weak def}). Thus,  $u (t, x)$ defined by (\ref{linear solution}) is a weak solution of the linear $G$-stochastic heat equation (\ref{linear equation}). 
  
  Finally, for a fixed initial condition $u_0$, suppose that $u$ and $v$ are two weak
  solutions of   $G$-stochastic heat equation (\ref{linear equation}).
  For any given $T>0$ and $\psi (t, x) \in C_c^{\infty} ([0, T] \times
  \mathbb{R})$, there exists a function $\phi (t, x) \in C_c^{\infty} ([0, T] \times
  \mathbb{R})$ satisfying the following partial differential equation:
  \begin{equation}
    \left\{\begin{array}{ll}
      & \frac{\partial}{\partial t} \phi (t, x) + \frac{\partial^2}{\partial
      x^2} \phi (t, x) = \psi (t, x),\ t > 0, x \in \mathbb{R},\\
      & \phi (T, x) = 0,\ x \in \mathbb{R} .
    \end{array}\right.
  \end{equation}
  Substituting $\phi$, $u$ and $v$ into  Eq. (\ref{linear weak def}) implies
  \begin{equation*}
   \int_0^T \int_{\mathbb{R}} (u (t, x) - v (t, x)) \psi (t, x) d x d t =0.
  \end{equation*}
  Due to the arbitrariness of $T$ and $\psi$, it is easy to check that  
  \begin{equation*} 
  u (t, x) = v (t, x) \text{ for a. e. } (t, x) \in [0, \infty)\times\mathbb{R},
  \end{equation*}
  which completes the proof.
\end{proof}

\subsection{Nonlinear $G$-stochastic heat equations}

In this subsection, for each fixed constants $T > 0, L > 0,$ we consider the following
nonlinear stochastic heat equation driven by the multiplicative space-time
$G$-white noise:
\begin{equation}
  \left\{\begin{array}{ll}
    & \frac{\partial}{\partial t} u (t, x) = \frac{\partial^2}{\partial x^2}
    u (t, x) +b(u(t,x))+ a (u (t, x)) \dot{\mathbf{W}} (t, x),\ 0 < t \leq T, 0 \leq x
    \leq L,\\
    & \frac{\partial}{\partial x} u (t, 0) = \frac{\partial}{\partial x} u
    (t, L) = 0,\ 0 < t \leq T,\\
    & u (0, x) = u_0 (x),\ 0 \leq x \leq L,
  \end{array}\right. \label{nonlinear eq}
\end{equation}
where  $u_0 : [0,L] \rightarrow \mathbb{R}$ is  a nonrandom, bounded and Borel measurable function, and  
 $a(x), b(x) \in C_{Lip} (\mathbb{R})$ are Lipschitz continuous functions, that is, there exists a constant $C$ such that 
\begin{equation*}
|a(x)-a(y)|+ |b(x)-b(y)|\leq C|x-y|, \text{ for all } x,y\in\mathbb{R}.
\end{equation*}

We first define the mild solution for  (\ref{nonlinear eq}) and subsequently demonstrate that the mild solution is also a weak solution.

\begin{definition}
  A spatio-temporal random field $\{u (t, x) : (t, x) \in [0, T] \times [0, L] \}$
  is said to be a {mild solution} of the nonlinear $G$-stochastic heat equation (\ref{nonlinear eq}) if it satisfies the
  following conditions:
  \begin{itemize}
    \item[(i)]  $(u (t, x))_{0 \leq t \leq T, 0 \leq x \leq L} \in \mathbf{S}_G^2
    ([0, T] \times [0, L])$;
    
    \item[(ii)] For $0<t\leq T$, $0\leq x\leq L$, 
    \begin{align}
         u  (t, x) = &\int_0^L u_0 (y) g (t, x, y) dy+ \int_0^t \int_0^L g (t - s, x, y) b (u(s, y))dyds \label{nonlinear mild def}\\
         &+\int_0^t \int_0^L g (t - s, x, y) a (u(s, y)) \mathbf{W} (ds, dy),\nonumber
    \end{align}
  \end{itemize}
  where $g (t, x, y)$ is the Green's  function corresponding to  Eq. (\ref{nonlinear eq}). 
\end{definition}

By Walsh \cite{Walsh} and Dalang et al. \cite{DK} , $g (t, x, y)$ is given by
\begin{eqnarray*}
  {g (t, x, y)} & = & \sum_{n = - \infty}^{\infty} [p (t, x - y - 2 {nL}) +
  p (t, x + y - 2 {nL})]\\
  & = & \frac{{1}}{\sqrt{4 \pi t}}  \sum_{n = - \infty}^{\infty} \left[ e^{-
  \frac{(x - y - 2 {nL})^2}{4 t}} + e^{- \frac{(x + y - 2 {nL})^2}{4
  t}} \right],\ 0 < t \leq T, 0 \leq x,y \leq L.
\end{eqnarray*}
Since for each fixed $0\leq x\leq L$, $g(t,x,y)\in L^2([0,T]\times[0,L])$, the integrals on the right-hand side of definition (\ref{nonlinear mild def}) make sense for $u\in \mathbf{S}_G^2([0,T]\times[0,L])$.

For each $\eta \in \mathbf{S}_G^2 ([0, T] \times [0, L])$, denote
\begin{equation}
(g \star \eta) (t, x) := \int_0^t \int_0^L g (t - s, x, y) \eta (s, y) \mathbf{W}
   (ds, dy),\ \text{for } 0 < t \leq T, 0 \leq x \leq L,
\end{equation}
and $(g \star \eta) (0, x) :=0$ for $0\leq x\leq L$. 
By Lemma \ref{Ito ineq}, $(g \star
\eta) (t, x) \in \mathbf{L}_G^2 (\mathcal{F}_t)$ for  $0 \leq t \leq T, 0 \leq x \leq L$.
Based on Lemma \ref{pW-est} and the properties of $g$, we obtain the following estimates and integrability for $g \star \eta$.

\begin{lemma}\label{geta}
  For each $\eta  \in \mathbf{S}_G^2 ([0, T] \times [0,
  L])$, there exists a constant $C$,  depending on $\overline{\sigma}^2$, $T$ and $L$, such that for any $0 \leq s, t \leq T$, $0\leq x, z \leq L,$
  \begin{equation} \label{g-eta-x}
  \hat{ \mathbb{E}} [| (g \star \eta) (t, x) - (g \star \eta) (t, z) |^2]
     \leq C {| x - z |}, 
  \end{equation}
  \begin{equation} \label{g-eta-t}
  \hat{ \mathbb{E}} [| (g \star \eta) (t, x) - (g \star \eta) (s, x) |^2]
     \leq C | t - s |^{1 / 2} . 
  \end{equation}
\end{lemma}
\begin{proof}
For each $t>0$, $x,z\in[0,L]$,
\begin{align*}
&\hat{\mathbb{E}} [|(g \star \eta) (t, x) - (g \star \eta) (t, z) |^2]\\
&=\hat{\mathbb{E}} \left[\left|\int_0^t\int_0^L\eta(s,y)(g(t-s,x,y)-g(t-s,z,y))\mathbf{W}(ds,dy)\right|^2\right]\\
&\leq \overline{\sigma}^2\|\eta\|^2_{\mathbf{S}^2}\int_0^t\int_0^L|g(t-s,x,y)-g(t-s,z,y)|^2dyds.
\end{align*}
By Lemma B.3.1 in Dalang and Sanz-Sol\'e \cite{DS}, we obtain that 
there exists a constant $C$ depending on $\overline{\sigma}^2$  such that
\begin{equation*}
\hat{\mathbb{E}} [|(g \star \eta) (t, x) - (g \star \eta) (t, z) |^2] \leq C {| x - z|},
\end{equation*}
which implies the estimate (\ref{g-eta-x}).

Similarly, for any $0\leq s\leq t\leq T$, we have
\begin{align*}
&\hat{\mathbb{E}} [| (g \star \eta) (t, x) - (g \star \eta) (s, x)|^2]\\
&\leq \overline{\sigma}^2\hat{\mathbb{E}} \left[\int_0^s\int_0^L|\eta(r,y)(g(t-r,x,y)-g(s-r,x,y))|^2dydr\right]\\
&\ \ \ +\overline{\sigma}^2\hat{\mathbb{E}} \left[\int_s^t\int_0^L\left|\eta(r,y)g(t-r,x,y)\right|^2dydr\right]\\
&\leq \overline{\sigma}^2\|\eta\|^2_{\mathbf{S}^2}\left(\int_0^s\int_0^L|g(t-r,x,y)-g(s-r,x,y)|^2dydr+\int_s^t\int_0^L\left|g(t-r,x,y)\right|^2dydr\right).
\end{align*}
Applying Lemma B.3.1 in Dalang and Sanz-Sol\'e \cite{DS} again,
there exists a constant $C$ depending on $T, L,\overline{\sigma}^2$  such that
\begin{equation*}
\hat{\mathbb{E}} [| (g \star \eta) (t, x) - (g \star \eta) (s, x) |^2] \leq C {|t - s|^{1/2}} .
\end{equation*}
\end{proof}

\begin{lemma}
  \label{gf}If  $\eta  \in \mathbf{S}_G^2 ([0, T] \times [0, L])$, then
  $g \star \eta \in \mathbf{S}_G^2 ([0, T] \times [0, L])$.
\end{lemma}

\begin{proof}
  For $n \geq 1$, define
  \begin{equation*} 
  g_n (t, x) := \sum_{i = 0}^{n - 1} \sum_{j = 0}^{n - 1} (g \star \eta)
     \left( \frac{i T}{n}, \frac{j L}{n} \right)  \textbf{1}_{[\frac{i T}{n}
     , \frac{(i + 1) T}{n})}(t)\textbf{1}_{[\frac{j L}{n} , \frac{(j + 1) L}{n}
     )}(x), \ \forall (t, x) \in [0, T] \times [0, L] . 
  \end{equation*}
  Then $g_n \in \mathbf{M}^{2, 0} ([0, T] \times [0, L])$ for $n\geq1$. By Lemma \ref{geta},
  there exists a constant $C$ depending on  $T$, $L$  and $\overline{\sigma}^2$ such that
  \begin{align*}
    &   \sup_{t \in [0, T]} \sup_{x \in [0, L]} \hat{ \mathbb{E}} [| (g
    \star \eta) (t, x) - g_n (t, x) |^2]\\
    &   = \sup_{t \in [0, T]} \sup_{x \in [0, L]} \hat{ \mathbb{E}}
    \left[ \sum_{i = 0}^{n - 1} \sum_{j = 0}^{n - 1} \left| (g \star \eta) (t, x)
    - (g \star \eta) \left( \frac{i T}{n}, \frac{j L}{n} \right) \right|^2 \textbf{1}_{[\frac{i T}{n}
     , \frac{(i + 1) T}{n})}(t)\textbf{1}_{[\frac{j L}{n} , \frac{(j + 1) L}{n}
     )}(x)\right]\\
    &   \leq \frac{ C}{\sqrt{n}} .
  \end{align*}
  Hence,
  \begin{equation*}
   \lim_{n \rightarrow \infty} \sup_{t \in [0, T]} \sup_{x \in [0, L]}
     \hat{ \mathbb{E}} [| (g \star \eta) (t, x) - g_n (t, x) |^2] = 0, 
  \end{equation*}
  which imples that $g \star \eta\in \mathbf{S}_G^2 ([0, T] \times [0, L]) .$
\end{proof}

Similar to Lemma \ref{geta} and Lemma \ref{gf}, we have the following statement.
\begin{lemma}
  \label{gb}For each  $\eta \in \mathbf{S}_G^2 ([0, T] \times [0, L])$, set
  \begin{equation*}
  \xi(t,x):=\int_0^t\int_0^L g(t-s,x,y)\eta(s,y)dyds,\ (t, x)\in(0,T]\times[0,L],
  \end{equation*}
and $\xi(0,x):=0$ for $x\in[0,L]$.  Then $\left\{\xi(t,x): (t, x)\in[0,T]\times[0,L]\right\}\in \mathbf{S}_G^2 ([0, T] \times [0, L])$.
\end{lemma}

\begin{theorem}\label{nonlinear-thm}
  Let $u_0 (x) $ be bounded and $a (x), b(x)$ be Lipschitz continuous functions. Then the nonlinear $G$-stochastic heat equation (\ref{nonlinear eq})  
  has a unique mild solution $\{u (t, x) : (t, x) \in [0, T] \times [0, L] \}$ such that $u  \in \mathbf{S}_G^2 ( [0, T] \times [0, L])$
  and $u (t, x)$ satisfies Eq. (\ref{nonlinear mild def}).
\end{theorem}

\begin{proof}
  (\textit{Existence}) For each $(t, x) \in (0, T] \times [0, L]$, $n\geq1$, define
  \begin{equation}
    \left\{\begin{array}{ll}
       u^0  (t, x) = &\int_0^L u_0 (y) g (t, x, y) dy,\\
       u^n  (t, x) = & u^0  (t, x) + \int_0^t \int_0^L g (t - s, x, y) b (u^{n
      - 1} (s, y)) dy ds\\
      &+ \int_0^t \int_0^L g (t - s, x, y) a (u^{n
      - 1} (s, y)) \mathbf{W} (ds, dy),
    \end{array}\right. \label{iteration}
  \end{equation}
  and $u^n(0,x):=u_0(x)$ for $0\leq x\leq L$, $n\in\mathbb{N}$. By Lemma \ref{gf} and Lemma \ref{gb}, the above iteration is
  well-defined and $\{u^n\}_{n=0}^\infty\subset\mathbf{S}_G^2 ( [0, T] \times [0, L])$. We then show that $u^n$ converges to a mild solution $u$ in
  $\mathbf{S}_G^2 ( [0, T] \times [0, L])$. For $n\geq1$,
  \begin{align*}
    &\hat{\mathbb{E}} [|u^{n + 1} (t, x) - u^n (t, x) |^2] \\
    & \leq 2 \hat{\mathbb{E}}
    \left[ \left| \int_0^t \int_0^L g (t - s, x, y) (b (u^n (s, y)) - b (u^{n
    - 1} (s, y))) dy ds\right|^2 \right]\\
    &\ \ \ +2\hat{\mathbb{E}}
    \left[ \left| \int_0^t \int_0^L g (t - s, x, y) (a (u^n (s, y)) - a (u^{n
    - 1} (s, y))) \mathbf{W} (ds, dy) \right|^2 \right]\\
    & \leq 2TL \hat{\mathbb{E}}
    \left[ \int_0^t \int_0^L g^2(t - s, x, y) \left| b (u^n (s, y)) - b (u^{n
    - 1} (s, y))\right|^2dy ds \right]\\
    &\ \ \ + 2\overline{\sigma}^2  \hat{\mathbb{E}} \left[ \int_0^t \int_0^L g^2 (t -
    s, x, y) |a (u^n (s, y)) - a (u^{n - 1} (s, y)) |^2 dyds \right]\\
    & \leq C  \int_0^t \int_0^L g^2 (t - s, x, y) 
    \hat{\mathbb{E}} [|u^n (s, y) - u^{n - 1} (s, y) |^2] dyds\\
    & \leq C    \int_0^t \sup_{y \in [0, L]} 
    \hat{\mathbb{E}} [|u^n (s, y) - u^{n - 1} (s, y) |^2]  \left[ \int_0^L g^2
    (t - s, x, y)  dy\right] ds,
  \end{align*}
  where $C$ is a constant depending on $a,b,T,L,\overline{\sigma}^2$.
  Since there exists a constant $C_{T,L}$ depending on $T$ and $L$ such that
    \begin{equation*}
  g(t,x,y)\leq \frac{C_{T,L}}{\sqrt{t}}e^{-\frac{|x-y|^2}{8t}},\ \forall 0<t\leq T,0\leq x,y\leq L,
    \end{equation*}
  then there exists a constant $C_1$ depending on $a,b,T,L,\overline{\sigma}^2$ satisfying
  \begin{equation}\label{ineq-H}
  \sup_{x \in [0, L]}  \hat{\mathbb{E}} [|u^{n + 1} (t, x) - u^n (t, x)
     |^2] \leq C_1  \int_0^t
     \sup_{y \in [0, L]}  \hat{\mathbb{E}} [|u^n (s, y) - u^{n - 1} (s, y)
     |^2]  (t - s)^{- \frac{1}{2}} ds.
 \end{equation}

  Denote $H_n (t) := \sup_{x \in [0, L]}  \hat{\mathbb{E}} [|u^{n + 1}
  (t, x) - u^n (t, x) |^2]$ for $0\leq t \leq T$, $n \in\mathbb{N}$. Equation (\ref{ineq-H}) reduces to 
  \begin{align}\label{ineq-H-2}
    H_n (t) & \leq C_1  \int_0^t
    H_{n - 1} (s)  (t - s)^{- \frac{1}{2}} ds \leq 4C_1^2  \int_0^t
    H_{n - 2} (s)  ds.
  \end{align} 
Applying induction, we have 
\begin{align}\label{ineq-H-3}
    H_n (t) & \leq  \frac{(4C_1^2)^{\lfloor {\frac{n}{2}}\rfloor}}{({\lfloor \frac{n}{2}\rfloor}-1)!}  \int_0^t
    H_{n - 2\lfloor \frac{n}{2}\rfloor} (s)(t-s)^{\lfloor \frac{n}{2}\rfloor-1}  ds,
  \end{align}
where $\lfloor {\frac{n}{2}}\rfloor$ denotes  the greatest integer less than or equal to $\frac{n}{2}$.
According to inequality (\ref{ineq-H-2}), if $H_n (t)$ is bounded
  on $[0, T]$,  $H_{n + 1} (t)$ is also bounded on $[0, T]$, $n \geq 0$. As $u^0$ is bounded, and 
  \begin{align*}
    H_0 (t) & \leq \sup_{x \in [0, L]} 2 \hat{\mathbb{E}} \left[ \left| \int_0^t
    \int_0^L g (t - s, x, y) b (u^0 (s, y)) dy ds  \right|^2\right.\\
    &\quad\quad\quad\quad\quad\quad\left.+ \left| \int_0^t
    \int_0^L g (t - s, x, y) a (u^0 (s, y)) \mathbf{W} (ds, dy) \right|^2
    \right]\\
    & \leq \sup_{x \in [0, L]} C_2  (\|u^0\|_{\mathbf{S}^2}^2+1) \int_0^t \int_0^L
    g^2 (t - s, x, y) dyds\\
    & \leq  C_3 (\|u^0\|_{\mathbf{S}^2}^2+1)\sqrt{ t},
  \end{align*}
   where $C_2, C_3$ are  constants depending on $a,b,T,L,\overline{\sigma}^2$, we know that  $H_n (t)$ is bounded on $[0, T]$ for $n \geq 0$.

   Denote $A:=\left|\sup_{0\leq t\leq T}H_0(t)\right|\vee\left|\sup_{0\leq t\leq T}H_1(t)\right|$.
    By  inequality (\ref{ineq-H-3}),   we deduce that
   \begin{align}\label{ineq-H-4}
   F^2_n(t):=\sup_{0\leq s\leq t} H_n (s) & \leq  \frac{A(4C_1^2)^{\lfloor {\frac{n}{2}}\rfloor}t^{\lfloor \frac{n}{2}\rfloor}}{{\lfloor \frac{n}{2}\rfloor}!},\ 0\leq t\leq T.
  \end{align} 
Then $\sum_{n=0}^\infty F_n(t)<\infty$, and we derive that 
  \begin{equation*} 
  \lim_{m, n \rightarrow \infty} \sup_{t \in [0, T]} \sup_{x \in [0, L]} 
     \hat{\mathbb{E}} [|u^m (t, x) - u^n (t, x) |^2] = 0.
  \end{equation*}
Thus, $\{u^n\}_{n=0}^\infty$ is a Cauchy sequence in $\mathbf{S}_G^2([0,T]\times [0,L])$.  For each $(t, x) \in [0, T] \times [0, L]$, set $u (t, x) =\mathbb{L}^2 -
  \lim_{n \rightarrow \infty} u^n  (t, x)$. Then $u \in \mathbf{S}_G^2
  ([0, T] \times [0, L]) .$

Letting $n\rightarrow \infty$ in the iteration (\ref{iteration}),  for $ 0<t\leq T$, $0\leq x\leq L$, we obtain that, in $\mathbf{L}_G^2(\mathcal{F}_t)$, 
  \begin{align*}
    u (t, x) = &\int_0^Lu_0(y)g(t,x,y)dy +\int_0^t \int_0^L g (t - s, x, y) b (u (s, y))  dy ds\\
    &+ \int_0^t \int_0^L g (t - s, x, y) a (u (s, y)) \mathbf{W} (ds, dy), 
  \end{align*}
 which proves the existence of the mild solution.

 (\textit{Uniqueness}) Suppose $v\in \mathbf{S}_G^2([0,T]\times [0,L]) $ is another mild solution of $G$-stochastic heat equation (\ref{nonlinear eq}) with the same initial condition $u_0$. Then  for $0<t\leq T$, $0\leq x\leq L$,
  \begin{align*}
    u (t, x)-v(t,x) = &\int_0^t \int_0^L g (t - s, x, y) [b(u (s, y))-b(v(s,y))] dyds\\
    &+ \int_0^t \int_0^L g (t - s, x, y) [a (u (s, y))-a(v(s,y))] \mathbf{W} (ds, dy).
  \end{align*}
Squaring and taking expectations on both sides yields
\begin{align*}
\hat{\mathbb{E}}[|u (t, x)-v(t,x)|^2]&\leq2\hat{\mathbb{E}}\left[\left|\int_0^t \int_0^L g (t - s, x, y) [b (u (s, y))-b(v(s,y))] dyds\right|^2\right]\\&\ \ \ +2\hat{\mathbb{E}}\left[\left|\int_0^t \int_0^L g (t - s, x, y) [a (u (s, y))-a(v(s,y))] \mathbf{W} (ds, dy)\right|^2\right]\\
&\leq C\int_0^t \int_0^L g^2 (t - s, x, y) \hat{\mathbb{E}}\left[|u (s, y)-v(s,y)|^2\right] dy ds,
\end{align*}
where $C$ is a constant depending on $a,b,T,L,\overline{\sigma}^2$.
As in the proof of existence,  we deduce that 
\begin{equation*} 
  \sup_{t \in [0, T]} \sup_{x \in [0, L]} 
     \hat{\mathbb{E}} [|u (t, x) - v (t, x) |^2] = 0, 
  \end{equation*}
which implies the uniqueness of the mild solution.
\end{proof}

Furthermore, we verify that the mild solution (\ref{nonlinear mild def}) also qualifies as a weak solution for  (\ref{nonlinear eq}).

\begin{definition}
  A spatio-temporal random field $\{u (t, x) : (t, x) \in [0, T] \times [0, L] \}$
  is said to be a {weak solution} of the $G$-stochastic heat equation (\ref{nonlinear eq})  if it satisfies the following conditions:
  \begin{itemize}
    \item[(i)]  $(u (t, x))_{0 \leq t \leq T, 0 \leq x \leq L} \in \mathbf{S}_G^2
    ([0, T] \times [0, L])$;
    
    \item[(ii)] For any $0<r\leq T$ and $\phi (t, x) \in C_c^{\infty} ([0, T] \times [0, L])$
    with $\frac{\partial}{\partial x} \phi (t, 0) = \frac{\partial}{\partial
    x} \phi (t, L) = 0$ for $0 \leq t \leq T$,
    \begin{align}
       & \int_0^L u (r, x) \phi (r, x) dx - \int_0^L u_0 (x) \phi (0, x) dx
      \nonumber\\
      & = \int_0^r \int_0^L u (t, x)  \left( \frac{\partial \phi (t,
      x)}{\partial t} + \frac{\partial^2 \phi (t, x)}{\partial x^2} \right)dxdt 
      + \int_0^r \int_0^L b(u(t,x))\phi (t, x) dxdt\label{nonlinear weak def}\\
       &\ \ \ + \int_0^r \int_0^L a(u(t,x))\phi (t, x) \mathbf{W} (dt, dx).\nonumber 
    \end{align}
  \end{itemize}
\end{definition}

\begin{theorem}
Let $u_0 (x) $ be bounded and $a(x), b(x)$ be Lipschitz continuous functions. Then the  mild solution, defined by (\ref{nonlinear mild def}), is also the weak solution of $G$-stochastic heat equation (\ref{nonlinear eq}).
\end{theorem}
\begin{proof}
Let $u \in \mathbf{S}_G^2([0, T] \times [0, L])$ be the mild solution of  (\ref{nonlinear eq}).
 It suffices to  prove that $u$ satisfies Eq. (\ref{nonlinear weak def}). For each $0<r\leq T$, 
  $\phi \in C_c^{\infty} ([0, T] \times [0, L])$ with $\frac{\partial}{\partial x} \phi (t, 0) = \frac{\partial}{\partial
    x} \phi (t, L) = 0$, $0 \leq t \leq T$, denote
  \begin{align}
    \mathcal{I}_1 := & \int_0^L u  (r, x) \phi (r, x) dx - \int_0^L u_0
    (x) \phi (0, x) dx, \label{substitute} \\
    \mathcal{I}_2 := & \int_0^r \int_0^L u (t, x)  \left(\frac{\partial
    \phi (t, x)}{\partial t} + \frac{\partial^2 \phi (t, x)}{\partial x^2}\right)
    dxdt, \label{I2} \\
    \mathcal{I}_3 := & \int_0^r \int_0^L b(u(t,x))\phi (t, x) dxdt+ \int_0^r \int_0^L a (u (t, x)) \phi (t, x) \mathbf{W} (dt, dx).\label{I3} 
  \end{align}
Then Eq. (\ref{nonlinear weak def}) is equal to $\mathcal{I}_1-\mathcal{I}_2-\mathcal{I}_3=0$.  Substituting Eq. (\ref{nonlinear mild def}) into $\mathcal{I}_1$ and $\mathcal{I}_2$ implies
  \begin{align*}
    \mathcal{I}_1 = & \int_0^L \int_0^L u_0 (y) g (r, x, y) \phi (r, x) dydx -
    \int_0^L u_0 (x) \phi (0, x) dx\\
    & + \int_0^L \left( \int_0^r \int_0^L g (r - s, x, y) b (u(s,
    y))dyds \right) \phi (r, x) dx\\
    & + \int_0^L \left( \int_0^r \int_0^L g (r - s, x, y) a (u(s,
    y)) \mathbf{W} (ds, dy) \right) \phi (r, x) dx\\
    \mathcal{I}_2 = & \int_0^r \int_0^L \left( \int_0^L u_0 (y) g (t, x, y) dy
    \right)  \left( \frac{\partial \phi (t, x)}{\partial t} + \frac{\partial^2
    \phi (t, x)}{\partial x^2} \right) dxdt\\
    & + \int_0^r \int_0^L \left( \int_0^t \int_0^L g (t - s, x, y) b (u(s, y)) dyds \right)  \left( \frac{\partial \phi (t,
    x)}{\partial t} + \frac{\partial^2 \phi (t, x)}{\partial x^2} \right)
    dxdt\\
    & + \int_0^r \int_0^L \left( \int_0^t \int_0^L g (t - s, x, y) a (u(s, y)) \mathbf{W} (ds, dy) \right)  \left( \frac{\partial \phi (t,
    x)}{\partial t} + \frac{\partial^2 \phi (t, x)}{\partial x^2} \right)
    dxdt.
  \end{align*}
  Since
  \begin{equation*}
   \int_0^L \int_0^L u_0 (y) g (0, x,y) \phi (0, x) dydx =
     \int_{0}^L u_0 (y) \phi (0, y) dy,
  \end{equation*}
  we have
  \begin{eqnarray*}
    &  & \int_0^L \int_0^r \int_0^L u_0 (y) g (t, x, y)  \left(
    \frac{\partial \phi (t, x)}{\partial t} + \frac{\partial^2 \phi (t,
    x)}{\partial x^2} \right) dxdtdy\\
    &  & = \int_0^L \int_0^L u_0 (y) g (r, x, y) \phi (r, x) dy d x -
    \int_0^L u_0 (y) \phi (0, y) dy.
  \end{eqnarray*}
  Substituting the above equality into $\mathcal{I}_1 -\mathcal{I}_2$ and applying Theorem \ref{Fubini}, we obtain
  \begin{align*}
    &\mathcal{I}_1 -\mathcal{I}_2 \\
    = & \int_0^r \int_0^L \int_0^L g (r - s, x,
    y) b (u(s, y)) \phi (r, x) dx dyds\\
    & - \int_0^r \int_0^L \int_s^r \int_0^L g (t - s, x, y) b (u(s,
    y))  \left( \frac{\partial \phi (t, x)}{\partial t} + \frac{\partial^2
    \phi (t, x)}{\partial x^2} \right) dxdt dy ds\\
    & +\int_0^r \int_0^L \int_0^L g (r - s, x,
    y) a (u(s, y)) \phi (r, x) dx \mathbf{W} (ds, dy)\\
    & - \int_0^r \int_0^L \int_s^r \int_0^L g (t - s, x, y) a (u(s,
    y))  \left( \frac{\partial \phi (t, x)}{\partial t} + \frac{\partial^2
    \phi (t, x)}{\partial x^2} \right) dxdt \mathbf{W} (ds, dy).
  \end{align*}
  
  As
  \begin{equation} \int_s^r \int_0^L g (t - s, x, y)  \left( \frac{\partial \phi (t,
     x)}{\partial t} + \frac{\partial^2 \phi (t, x)}{\partial x^2} \right)
     dxdt = \int_0^L g (r - s, x, y) \phi (r, x) dx - \phi (s, y), 
  \end{equation}
  we get  $\mathcal{I}_1 -\mathcal{I}_2 = \mathcal{I}_3$, 
  which completes the proof.
\end{proof}

At the end of this section, we  present the moment estimates for the  solution. Utilizing Theorem 4.3 in Sanz-Sol\'e  and Sarr\`a \cite{Sanz} and Lemma \ref{geta}, we derive the following estimates. The proof is straightforward and  omitted here.
\begin{proposition}
  Assume that $u_0(x):[0, L]\rightarrow \mathbb{R} $ is a bounded and $\alpha$-H\"older continuous function for some $\alpha>0$, and $a(x), b(x)$ are Lipschitz continuous functions. Let $\{u (t,
  x) : (t, x) \in [0, T] \times [0, L] \}$ be the mild solution of  $G$-stochastic heat equation (\ref{nonlinear eq}). 
  Then  there exists a constant $C$ such that for any $ s,t \in [0, T]$, $x, z \in [0, L],$
  \begin{itemize}
    \item[(i)] $\hat{\mathbb{E}} [|u (t, x) - u (t, z) |^2] \leq C {| x - z|^{1\wedge(2\alpha)}} ;$
    
    \item[(ii)] $\hat{\mathbb{E}} [|u (t , x) - u (s, x) |^2] \leq C |
    t - s |^{(1 / 2)\wedge\alpha} .$
  \end{itemize}
\end{proposition}

\section{Applications to problems with probability uncertainty}

In many real-world problems, due to  inherent complexities and model uncertainty, space-time $G$-white noise may provide a more appropriate approximation to  the  random noise, making   $G$-stochastic heat equations better suited to describe  real random systems. In this section, we present some examples to illustrate the potential applications of  $G$-stochastic heat equations.

\begin{example}[Random motion of a polymer chain in a liquid]\label{e1}
Consider a polymer chain in a liquid. The continuum limit of the random motion of  particles satisfies the following stochastic heat equation (see Hairer \cite{Hai}),
\begin{equation}\label{Eg1}
  \left\{\begin{array}{ll}
    & \frac{\partial}{\partial t} u (t, x) =\frac{\partial^2}{\partial x^2}
    u (t, x) +b(u(t,x))+  \xi(t, x),\  t>0 , 0 \leq x \leq 1,\\
    & \frac{\partial}{\partial x}u (t, 0) =  \frac{\partial}{\partial x}u(t, 1) = 0,\  t >0,
  \end{array}\right.
\end{equation}
where $b(u(t,x))$ denotes the deterministic external force and $\xi(t,x)$ describes the random kicks from the surrounding liquid. 
Assuming the random kicks occur independently at a high rate and their distribution is identical  and invariant, $\xi(t,x)$ is  modeled as a classical space-time white noise using the central limit theorem for independent and identically distributed random variables.

However, when the surrounding liquid is not constant (e.g., the temperature   fluctuates within an interval),  the distribution of kicks also fluctuates in a set. In this situation, by Peng's central limit theorem under sublinear expectation \cite{P2010}, $G$-normally distributed random variables may be a more suitable approximation for modeling $\xi(t,x)$. Thus, $\xi(t,x)$ can be replaced by  $\dot{\mathbf{W}}(t,x)$, transforming Eq. (\ref{Eg1}) into the $G$-stochastic heat equation (\ref{nonlinear eq}) with the coefficient $a\equiv1$.
\end{example}

\begin{example}[Heat density in a random medium]

In \cite{Kh14}, Khoshnevisan introduced the following stochastic heat equation with Dirichlet boundary conditions to describe the heat flow in a rod of unit length, 
\begin{equation}\label{Eg2}
  \left\{\begin{array}{ll}
    & \frac{\partial}{\partial t} u (t, x) = \nu\frac{\partial^2}{\partial x^2}
    u (t, x) + a (u (t, x)) \xi(t, x),\  t>0 , 0 \leq x \leq 1,\\
    & u (t, 0) =  u(t, 1) = 0,\  t >0,
  \end{array}\right.
\end{equation}
where $\nu$ is a positive constant,  $a(x)$ is a smooth function characterizing the feedback mechanism, and the space-time white noise $\xi(t,x)$ describes the external forcing density. The solution $u(t,x)$ represents the heat density at $(t,x)$.

In real experiments, it is challenging to maintain a constant  external heating or cooling source. When the external source fluctuates within a set,   $\dot{\mathbf{W}}(t,x)$ can be used to replace $\xi(t,x)$ to account for the model uncertainty.  Consequently, Eq. (\ref{Eg2}) becomes a $G$-stochastic heat equation. 

It is worth mentioning that although the boundary conditions of Eqs. (\ref{Eg2}) and (\ref{nonlinear eq}) differ, by updating the Green's function in the solution,  the well-posedness of the solution in Theorem \ref{nonlinear-thm} can be extended to the $G$-stochastic heat equation with Dirichlet boundary conditions.
\end{example}

\begin{example}[Propagation of electric potentials in neurons]
Assume nerve cells are long and thin cylinders from $0$ to $L$. Let $V(t,x)$ represent the electrical potential at time $t$ and position $x$. To describe the propagation of $V(t,x)$, Walsh \cite{Walsh} introduced the following stochastic heat equation
\begin{equation}\label{Eg3}
\frac{\partial}{\partial t} V (t, x) = \frac{\partial^2}{\partial x^2}
    V (t, x) - V(t,x)+ f(V(t,x))\xi(t, x),\  t>0, 0 \leq x \leq L, 
\end{equation}
where $f(V(t,x))$ denotes the neuron's response  to a current impulse, and  the  space-time white noise $\xi(t,x)$ describes the random impulse.
In this model, the impulses are assumed to be small, numerous, and independent at different position $x$ and time $t$.

Similar to Example \ref{e1}, when 
the  distribution of impulses fluctuates within a set, the $G$-stochastic heat equation (\ref{nonlinear eq}) with the multiplicative noise $\dot{\mathbf{W}}(t,x)$ may be a more suitable model to describe the propagation of electric potentials in neurons.

\end{example}

\begin{example}[Parabolic Anderson models with probability uncertainty]
As a special case of the stochastic heat equation, the parabolic Anderson model can be studied within  the sublinear expectation framework. The results in this paper can be applied to such models, describing by the following equation with suitable initial and boundary conditions,
\begin{equation}\label{Eg4}
\frac{\partial}{\partial t} u (t, x) = \frac{\partial^2}{\partial x^2}
    u (t, x) + u(t,x)\dot{\mathbf{W}}(t, x),\  t>0, x\in D,
\end{equation}
where $D\subset\mathbb{R}$ is a domain with
smooth boundary.

\end{example}

%
%


%

\section*{Acknowledgements}
This work was partially supported by National Natural Science Foundation of China (No. 12401179), Natural Science Foundation of Shandong Province (No.  ZR2022QA012), and the Fundamental Research Funds for the Central Universities.

\end{document}